\documentclass[12pt]{article}
\usepackage{epsfig}
\usepackage{amssymb}
\usepackage{amsmath}

\newtheorem{thm}{Theorem}[section] 
\newtheorem{lemma}[thm]{Lemma} 
\newtheorem{Def}[thm]{Definition} 
 
\newtheorem{prop}[thm]{Proposition} 
\newtheorem{rmk}[thm]{Remark}

\newcommand{\ov}{\overline}
\renewcommand{\dfrac}{\displaystyle\frac}
\newcommand{\D}{{\mathbb D}}

\newcommand{\C}{{\mathbb C}}
\renewcommand{\mod}{\mbox{mod\,}}
\newcommand{\re}{\mbox{Re}\,}
\newcommand{\im}{\mbox{Im}\,}
\newcommand{\proof}{\noindent {\bf Proof} \hspace{0.2in}} 
\newcommand{\qed}{\hfill\mbox{\raggedright\rule{.07in}{.1in}}
  \vspace{1ex}} 
\newcommand{\dps}{\displaystyle}
\newcommand{\Section}[1]{\section{#1} \setcounter{equation}{0}}

\newcommand{\arraystart}{\renewcommand{\arraystretch}{1.6}}
\newcommand{\arrayfinish}{\renewcommand{\arraystretch}{1.2}}
 
\title{Equivariant versal unfoldings for
linear retarded functional differential 
equations}
 
\author{Pietro-Luciano Buono\\Centre de Recherches
  Math\'ematiques\\Universit\'e de Montr\'eal\\C.P. 6128, 
Succ. Centre-Ville\\ 
     Montr\'eal, QC H3C 3J7\\CANADA  \and Victor G. LeBlanc\\Department of
  Mathematics and Statistics\\University of Ottawa\\Ottawa, 
ON K1N 6N5\\CANADA}
\date{April 11, 2003}

\begin{document}

\maketitle

\begin{abstract}
We continue our investigation of versality for parametrized families of 
linear retarded functional differential equations (RFDEs) projected onto
finite-dimensional invariant manifolds. In this paper, we consider RFDEs 
equivariant with respect to the action of a compact Lie group. In a previous 
paper (Buono and LeBlanc, to appear in {\em J. Diff. Eqs.}), we have studied 
this question in the general case (i.e. no {\em a priori} restrictions on 
the RFDE). When studying the question of versality in the equivariant 
context, it is natural to want to restrict the range of possible unfoldings 
to include only those which share the same symmetries as the original RFDE, 
and so our previous results do not immediately apply. In this paper, we 
show that with appropriate projections, our previous results on versal 
unfoldings of linear RFDEs can be adapted to the case of linear equivariant 
RFDEs. We illustrate our theory by studying the linear equivariant 
unfoldings at double Hopf bifurcation points in a $\D_3$-equivariant network 
of coupled identical neurons modeled by delay-differential equations due 
to delays in the internal dynamics and coupling.
\end{abstract}

\Section{Introduction}

Symmetry plays an important role in the description of many physical
systems.  Whether the symmetry is an intrinsic part of the physical system, or
whether it is merely a modeling assumption, 
the mathematical model for any such system must reflect this
symmetry.  
When the mathematical model is a differential equation, 
this is achieved by imposing {\em equivariance} of the differential
equation.  That is, if $\Gamma$ is a group of transformations
which represents the
symmetries of the physical system, 
then the differential equation model must commute with an
action of $\Gamma$.  
An important consequence of $\Gamma$-equivariance is that
any solution of the differential
equation is mapped onto other solutions by the 
elements of $\Gamma$.  Thus, the algebraic structure of the group
$\Gamma$ yields important information about the geometrical
structure of solutions of the differential equation in phase space.

The theory of equivariant dynamical systems is by now
well-developed and understood (see, for example 
\cite{GS83} and \cite{GSS88}).  $\Gamma$-equivariance imposes
strong algebraic restrictions on the functional form of the 
differential equation.  This has important implications, for example, in the
study of local stability and bifurcation of equilibrium points.  In
particular, since the linearization of the vector field at the
equilibrium point must commute with $\Gamma$, then
this linearization usually has eigenvalues of high
geometric multiplicity.

In this paper, we will be interested in systems which are modeled by
parametrized families of retarded functional differential equations 
\begin{equation}
\dot{z}(t)={\cal L}(\alpha)\,z_t+F(z_t,\alpha),\quad\alpha\in {\mathbb C}^p,
\label{eq1}
\end{equation}
where ${\cal L}(\alpha)$ is a parametrized family of bounded linear
functional operators and $(u,\alpha)\mapsto F(u,\alpha)$ is nonlinear 
and smooth with $F(0,\alpha)=0$ and $D_{u}F(0,\alpha)=0$ for all $\alpha$, 
such that (\ref{eq1}) is equivariant
under the action of a compact Lie group of symmetries $\Gamma$ for
all $\alpha$.
Examples of such systems can be found in the study of neural networks
modeled by symmetrically coupled cell systems with delays in the 
coupling~(\cite{KVW98,NCW,Wu98,WFH99}), in the study of semiconductor lasers
subject to optical feedback~\cite{KTG00,VK00} or symmetrically coupled 
arrays of lasers~\cite{Heiletal01}.  
Specifically, we are interested in the following question regarding 
(\ref{eq1}): consider the equilibrium $x=0$, $\alpha=0$ of (\ref{eq1}), 
and let $V$ be a finite-dimensional invariant manifold for the dynamics 
of (\ref{eq1}) which passes through this equilibrium point.  Let
\begin{equation}
\dot{x}=B(\alpha)x+G(x,\alpha)
\label{eq2}
\end{equation}
be the parametrized family of ordinary differential equations which
gives the dynamics of (\ref{eq1}) restricted to $V$, where $B(\alpha)$ is
a parametrized family of matrices, and $G$ is the nonlinearity.
Under what conditions on ${\cal L}(\alpha)$ in (\ref{eq1}) can we be assured that
$B(\alpha)$ is a versal unfolding of the
matrix $B(0)$ {\em within the space of $\Gamma$-equivariant matrices}?
Recall that, loosely speaking, the family $B(\alpha)$ is a versal
unfolding of $B(0)$ if it is, up to similarity transformations and
changes of parameters, the
most general $C^{\infty}$ perturbation of $B(0)$ (a precise definition is
given in \cite{BL02} in the absence of $\Gamma$-equivariance, and in
Section \ref{sec3} of this paper in the $\Gamma$-equivariant case).
The above question is especially important in the case
where (\ref{eq2}) represents a center manifold reduction of
(\ref{eq1}).  If $B(\alpha)$ is not versal, then this could lead to
restrictions on the bifurcation behavior of (\ref{eq1}) near the equilibrium.
It is clear that answering this question does not depend on the
nonlinearity $F$ in (\ref{eq1}); therefore, we restrict our attention
to the class of parametrized families of 
$\Gamma$-equivariant linear retarded functional differential equations.

This problem was studied in \cite{BL02} in the
absence of $\Gamma$-equivariance.  In particular, 
Theorem 7.4 of that paper gives a method for constructing a
parametrized family ${\cal L}(\alpha)$ of bounded linear functional
operators such that the family of matrices $B(\alpha)$ in (\ref{eq2})
is a versal unfolding of the matrix $B(0)$.  Our approach in this
paper will be to construct a suitable equivariant projection of the family
${\cal L}(\alpha)$ of Theorem 7.4 of \cite{BL02}, and to show
that the resulting associated $\Gamma$-equivariant family of matrices
$B(\alpha)$ in (\ref{eq2}) is a versal unfolding of $B(0)$ within the 
space of all $\Gamma$-equivariant matrices.

The paper is organized as follows. In Section 2 we review some 
basic theory of linear retarded functional equations while Section 3
contains the results concerning versal and mini-versal unfoldings
of equivariant matrices. In Section 4, our main result about the 
construction of equivariant unfoldings for linear retarded functional
differential equations is presented. Finally, in Section 5 a detailed
computation of the equivariant mini-versal unfolding of a $\D_3$-equivariant
linear delay-differential equation at a double Hopf bifurcation is discussed
in both cases of simple and double imaginary eigenvalues.

It is assumed that the reader is familiar with the terminology and
results of \cite{BL02}, \cite{FMH}, \cite{FMR}, \cite{FM96} and
\cite{HalVL}.

\Section{Reduction of equivariant linear RFDEs}

Let $C_n=C([-\tau,0],{\mathbb C}^n)$ be the Banach space of continuous 
functions from the interval
$[-\tau,0]$, into ${\mathbb C}^n$ ($\tau>0$) endowed with the norm of 
uniform convergence.
We will first consider the linear homogeneous RFDE
\begin{equation}
\dot{z}(t)={\cal L}_0(z_t),
\label{fde}
\end{equation}
where ${\cal L}_0$ is a bounded linear operator from $C_n$ into 
${\mathbb C}^n$.  We write
\[
{\cal L}_0(\varphi)=\int_{-\tau}^{0}\,d\eta(\theta)\varphi(\theta),
\]
where $\eta$ is an $n\times n$ matrix-valued function of bounded variation 
defined on $[-\tau,0]$.  

We will suppose that (\ref{fde}) has certain
symmetry properties.  Specifically, we
let $\Gamma$ be a compact group of transformations acting linearly on
${\mathbb C}^n$.  We say that
(\ref{fde}) is {\em $\Gamma$-equivariant} if
\begin{equation}
\gamma\cdot d\eta(\theta)=d\eta(\theta)\cdot\gamma,\,\,\,\forall\,
\gamma\in\Gamma, \,\theta\in [-\tau,0].
\label{equiv_cond}
\end{equation}
A consequence of this definition is that if $x(t)$ is a solution to
(\ref{fde}), then so is $\gamma\cdot x(t)$, for all $\gamma\in\Gamma$.

Suppose that $\Lambda\subset {\mathbb C}$ is a non-empty finite set of
eigenvalues of the infinitesimal generator $A_0$ for the semi-flow
of (\ref{fde}), i.e. $\lambda\in\Lambda$ if and only if
$\lambda$ satisfies the
characteristic equation
\begin{equation}
\mbox{\rm det}\,\Delta(\lambda)=0,\;\;\;\;\;\;\;\;
\mbox{\rm where}\;\;\Delta(\lambda)=\lambda\,I_n-\int_{-\tau}^0\,
d\eta(\theta)e^{\lambda\theta},
\label{char_eq}
\end{equation}
where $I_n$ is the $n\times n$ identity matrix.
Using adjoint theory, it is known that we can write
\begin{equation}
C_n=P\oplus Q
\label{splitcn}
\end{equation}
where the generalized eigenspace $P$ corresponding to $\Lambda$ and
the complementary subspace
$Q$ are both invariant under the semiflow of (\ref{fde}), and invariant
under $A_0$.

Define $C_n^*=C([0,\tau],{\mathbb C}^{n*})$, where ${\mathbb C}^{n*}$ is 
the $n$-dimensional space of row vectors.  We have the adjoint bilinear 
form on $C_n^*\times C_n$:
\begin{equation}
(\psi,\varphi)_n=\psi(0)\varphi(0)-\int_{-\tau}^0\int_0^{\theta}\psi(\xi
-\theta)d\eta(\theta)\varphi(\xi)d\xi.
\label{bilin_form}
\end{equation}
We let $\Phi=(\varphi_1,\ldots,\varphi_c)$ be a basis for $P$, and 
$\Psi=\mbox{\rm col}(\psi_1,\ldots,\psi_c)$ be a basis for the dual 
space $P^*$ in $C_n^*$, chosen so that $(\Psi,\Phi)_n$ is the $c\times c$ 
identity matrix, $I_c$.  In this case, we have $Q=\{\varphi\in C_n\,:\,
(\Psi,\varphi)_n=0\}$.  
It follows that $\Phi'=\Phi B$ where $B$ is 
a $c\times c$ constant matrix.
The spectrum of $B$ coincides with $\Lambda$.
Using the decomposition (\ref{splitcn}), any $z\in C_n$ can be written as
$z=\Phi\,x+y$, where $x\in {\mathbb C}^c$ and $y\in Q$ is
a $C^1$ function.  The dynamics of (\ref{fde}) on $P$ are then given
by
the $c$-dimensional linear ordinary differential equation 
\begin{equation}
\dot{x}=Bx.
\label{red_ode}
\end{equation}

The important point, as we will now show, is that $\Gamma$-equivariance of (\ref{fde})
forces $\Gamma$-invariance of the splitting (\ref{splitcn}) and 
$\Gamma$-equivariance of
the reduced ordinary differential equation (\ref{red_ode}).

\begin{prop}
The spaces $P$, $Q$ and $P^*$ defined above are $\Gamma$ invariant, 
and the $c\times c$ matrix $B$ in (\ref{red_ode})
commutes
with a representation $G:\Gamma\longrightarrow GL(c,{\mathbb C})$ of the group $\Gamma$.
\label{basic_sym_prop}
\end{prop}
\proof
We start by noting the following identity which follows trivially from
$\Gamma$-equivariance of (\ref{fde}), and from (\ref{bilin_form}):
\begin{equation}
(\Psi,\gamma\cdot\Phi)_n=(\Psi\cdot\gamma,\Phi)_n\,\,\,\,\,\forall\,\gamma\in\Gamma.
\label{bilin_equiv}
\end{equation}

The functions $\varphi_1,\ldots,\varphi_c$ which span $P$ are
generalized eigenfunctions of the operator $A$ defined by
$A\,\varphi=\varphi'$ whose domain is
\[
{\cal D}(A)=\left\{\,\varphi\in C_n\,:\,\varphi'\in C_n\,\,\mbox{\rm and}\,\,
\varphi'(0)=\int_{-\tau}^0\,d\eta(\theta)\,\varphi(\theta)\,\right\}.
\]
Let $\gamma\in\Gamma$, $\varphi\in P$, and consider
$\tilde{\varphi}=\gamma\cdot \varphi$.  Then obviously
$\tilde{\varphi}'\in C_n$ since $\varphi'\in C_n$.  Moreover, from
(\ref{equiv_cond}) and the fact that $\varphi\in {\cal D}(A)$, we get
\[
\tilde{\varphi}'(0)=\gamma\cdot\varphi'(0)=\gamma\,\int_{-\tau}^0\,d\eta(\theta)\,\varphi(\theta)=
\int_{-\tau}^0\,d\eta(\theta)\,\gamma\cdot\varphi(\theta)=
\int_{-\tau}^0\,d\eta(\theta)\,\tilde{\varphi}(\theta),
\]
so $\tilde{\varphi}\in {\cal D}(A)$.  Now, there exists an integer $k$
and a $\lambda\in\Lambda$
such that $(A-\lambda I)^k\varphi=0$.  Since $A\,(\gamma\cdot\varphi) =
(\gamma\cdot\varphi)'=\gamma\cdot\varphi'=\gamma\cdot A\,\varphi$, it
is easy to see that $(A-\lambda I)^k\tilde{\varphi}=0$, and we
conclude that $\tilde{\varphi}\in P$.  Thus, $P$ is
$\Gamma$-invariant.  A similar argument shows that $P^*$ is also
$\Gamma$-invariant.  This implies that for all $\gamma\in\Gamma$,
there
exist $c\times c$ matrices $G(\gamma)$ and $H(\gamma)$ such that
\begin{equation}
\gamma\cdot\Phi=\Phi\cdot G(\gamma)\,\,\,\,\mbox{\rm and}\,\,\,\,
\Psi\cdot\gamma=H(\gamma)\cdot\Psi.
\label{c_act}
\end{equation}
In fact, it follows
from (\ref{bilin_equiv}) and the fact that $(\Psi,\Phi)_n=I_c$ that
\[
\begin{array}{l}
H(\gamma)=H(\gamma)(\Psi,\Phi)_n=(H(\gamma)\cdot\Psi,\Phi)_n=(\Psi\cdot\gamma,\Phi)_n=
(\Psi,\gamma\cdot\Phi)_n=\\[0.15in]
(\Psi,\Phi\cdot G(\gamma))_n=(\Psi,\Phi)_nG(\gamma)=G(\gamma).
\end{array}
\]
Note that
\[
\begin{array}{l}
G(\sigma\gamma)=(\Psi,\Phi)_nG(\sigma\gamma)=(\Psi,\Phi\cdot
G(\sigma\gamma))_n=(\Psi,\sigma\gamma\cdot\Phi)_n=\\[0.15in]
(\Psi\cdot\sigma,\gamma\cdot\Phi)_n=(G(\sigma)\cdot\Psi,\Phi\cdot
G(\gamma))_n = G(\sigma)(\Psi,\Phi)_nG(\gamma)=G(\sigma)G(\gamma),
\end{array}
\]
so that the mapping
\[
\begin{array}{clc}
\Gamma&\longrightarrow&GL(c,{\mathbb C})\\
\gamma&\longmapsto&G(\gamma)
\end{array}
\]
is a representation of the group $\Gamma$ in the space of invertible $c\times
c$ matrices.

Now, $\varphi\in Q$ if and only if $(\Psi,\varphi)_n=0$.  Then using
(\ref{bilin_equiv}), one gets that
\[
(\Psi,\gamma\cdot\varphi)_n=(\Psi\cdot\gamma,\varphi)_n=
(G(\gamma)\cdot\Psi,\varphi)_n=G(\gamma)(\Psi,\varphi)_n=0,
\]
so $Q$ is $\Gamma$-invariant.

Finally, it follows from the fact that $B=(\Psi,\Phi')_n$ that
\begin{equation}
B\cdot G(\gamma)=G(\gamma)\cdot B\,\,\,\,\,\,\,\,\forall\,\gamma\in\Gamma.
\label{BGGB}
\end{equation}
\hfill
\qed

\Section{Equivariant versal unfoldings of matrices\label{sec3}} 

Let $B\in \mbox{\rm
  Mat}_{c\times c}$, where $\mbox{\rm Mat}_{c\times c}$ denotes the space of $c\times c$
matrices with entries in ${\mathbb C}$.  
Recall that a $p$-{\em parameter unfolding} of $B$ is a $C^{\infty}$
mapping
${\cal B}:{\mathbb C}^p\longrightarrow\mbox{\rm Mat}_{c\times c}$ such
that ${\cal B}(\alpha_0)=B$ for some $\alpha_0\in {\mathbb C}^p$.
We then say that a $p$-parameter unfolding
${\cal B}(\alpha)$ of $B$ is a {\em
  versal unfolding} of $B$ if, given any $q$-parameter unfolding 
${\cal A}(\beta)$
of $B$ (with ${\cal A}(\beta_0)=B$), there exists a $C^{\infty}$ mapping
\[
\phi:{\mathbb C}^q\longrightarrow {\mathbb C}^p
\] 
with $\phi(\beta_0)=\alpha_0$,
and a $q$-parameter unfolding $C(\beta)$ of the identity matrix such
that
\[
{\cal A}(\beta)=C(\beta){\cal
  B}(\phi(\beta))(C(\beta))^{-1},\,\,\,\,\,\,\,\,\,\forall\,
  \beta\,\,\mbox{\rm near}\,\,\beta_0.
\]
If ${\cal B}$ is a versal unfolding which 
depends on the least number of parameters, then ${\cal B}$ is called a
  {\em mini-versal unfolding}.
Thus, one may view a versal unfolding of $B$ as the most general
$C^{\infty}$ perturbation of $B$ up to similarity and change of
parameters.
One then has the following sufficient criterion \cite{Ar} for establishing the
versality of a given unfolding of $B$:
\begin{prop}
Let ${\cal B}(\alpha)$ be a $p$-parameter unfolding of $B\in\mbox{\rm
  Mat}_{c\times c}$.  Let $\Sigma$ denote the similarity orbit of $B$,
  i.e.
\[
\Sigma=\{\,CBC^{-1}\,:\,C\in GL(c,{\mathbb C})\,\}\subset\mbox{\rm
  Mat}_{c\times c}.
\]
If
\begin{equation}
\mbox{\rm Mat}_{c\times c}=T_B\Sigma+D_{\alpha}{\cal B}(\alpha_0)\cdot
{\mathbb C}^p,
\label{suff_cond_versal} 
\end{equation}
where $T_B\Sigma$ is the tangent space to $\Sigma$ at $B$, then
${\cal B}$ is a versal unfolding of $B$.  
If, in addition, 
the codimension of $T_B\Sigma$ in $\mbox{\rm Mat}_{c\times c}$ is equal to $p$,
then ${\cal B}$ is mini-versal.
\label{suff_cond_versal_prop}
\end{prop} 

In this paper, we are not interested in the most general
perturbation of $B$, since we are only considering matrices in
$\mbox{\rm Mat}_{c\times c}$ which commute with the group $\Gamma$.
In fact, we need only consider perturbations in the subspace
$\mbox{\rm Mat}^{\Gamma}_{c\times c}$ defined by
\[
\mbox{\rm Mat}^{\Gamma}_{c\times c}=\{\,A\in\mbox{\rm Mat}_{c\times
  c}\,:\,
G(\gamma)A=AG(\gamma),\,\,\,\forall\,\gamma\in\Gamma\,\},
\]
where $G$ is the representation of $\Gamma$ in $GL(c,{\mathbb C})$ 
defined in (\ref{c_act}).
We therefore adopt the following definition:
\begin{Def}
A $p$-parameter unfolding ${\cal B}(\alpha)$ of
$B\in\mbox{\rm Mat}^{\Gamma}_{c\times c}$ is said to be a
{\em $\Gamma$-unfolding} of $B$ if ${\cal B}(\alpha)\in\mbox{\rm Mat}^{\Gamma}_{c\times c}$
for all $\alpha$.  The $\Gamma$-unfolding is said to be {\em $\Gamma$-versal} if,
given any $q$-parameter $\Gamma$-unfolding 
${\cal A}(\beta)$
of $B$ (with ${\cal A}(\beta_0)=B$), there exists a $C^{\infty}$ mapping
\[
\phi:{\mathbb C}^q\longrightarrow {\mathbb C}^p
\] 
with $\phi(\beta_0)=\alpha_0$,
and a $q$-parameter $\Gamma$-unfolding $C(\beta)$ of the identity matrix such
that
\[
{\cal A}(\beta)=C(\beta){\cal
  B}(\phi(\beta))(C(\beta))^{-1},\,\,\,\,\,\,\,\,\,\forall\,
  \beta\,\,\mbox{\rm near}\,\,\beta_0.
\]
If ${\cal B}$ is a $\Gamma$-versal unfolding which 
depends on the least number of parameters, then ${\cal B}$ is called a
  {\em $\Gamma$-mini-versal unfolding}.
\label{gammaunfolddef}
\end{Def}
The following lemma, which is the $\Gamma$-equivariant version of
Proposition \ref{suff_cond_versal_prop}, is proved in Appendix A.
\begin{lemma}
Let ${\cal B}(\alpha)$ be a $p$-parameter $\Gamma$-unfolding of $B\in\mbox{\rm
  Mat}_{c\times c}^{\Gamma}$.  Let $\Sigma^{\Gamma}$ 
denote the $\Gamma$-similarity orbit of $B$,
  i.e.
\[
\Sigma^{\Gamma}=\{\,CBC^{-1}\,:\,C\in GL(c,{\mathbb C})\cap\mbox{\rm Mat}_{c\times c}^{\Gamma}\,\}\subset\mbox{\rm
  Mat}_{c\times c}^{\Gamma}.
\]
If
\begin{equation}
\mbox{\rm Mat}_{c\times c}^{\Gamma}=T_B\Sigma^{\Gamma}+D_{\alpha}{\cal B}(\alpha_0)\cdot
{\mathbb C}^p,
\label{suff_cond_versal_gamma} 
\end{equation}
where $T_B\Sigma^{\Gamma}$ is the tangent space to $\Sigma^{\Gamma}$ at $B$, then
${\cal B}$ is a $\Gamma$-versal unfolding of $B$.  
If, in addition, 
the codimension of $T_B\Sigma^{\Gamma}$ in $\mbox{\rm Mat}^{\Gamma}_{c\times c}$ is
equal to $p$, then ${\cal B}$ is $\Gamma$-mini-versal.
\label{suff_cond_versal_prop_gamma}
\end{lemma}

\begin{lemma}
There exists a projection $\pi_c^{\Gamma}:\mbox{\rm Mat}_{c\times
  c}\longrightarrow \mbox{\rm Mat}^{\Gamma}_{c\times c}$ such that
\begin{equation}
\pi_c^{\Gamma}(AM)=A\,\pi_c^{\Gamma}(M)\,\,\,\,\,\mbox{\rm and}
\,\,\,\,\,\pi_c^{\Gamma}(MA)=\pi_c^{\Gamma}(M)\,A,\,\,\,\,\,\forall\,
A\in\mbox{\rm Mat}^{\Gamma}_{c\times c},\,\,M\in\mbox{\rm Mat}_{c\times c}.
\label{morphism}
\end{equation}
\label{homlem}
\end{lemma}
\proof
For $M\in\mbox{\rm Mat}_{c\times c}$, define
\[
\pi_c^{\Gamma}(M)=\int_{\Gamma}\,G(\gamma)MG(\gamma^{-1})\,d\gamma
\]
where the integral denotes the Haar integral on $\Gamma$, which we
assume
has been normalized so that the measure of $\Gamma$ is 1 (recall that
we assume that $\Gamma$ is compact).
Since $G$ is a morphism and since the Haar integral is translation
invariant,
we have that for all $\sigma\in\Gamma$
\[
G(\sigma)\pi_c^{\Gamma}(M)(G(\sigma))^{-1}=
\int_{\Gamma}\,G(\sigma\gamma)MG((\sigma\gamma)^{-1})\,d\gamma=\pi_c^{\Gamma}(M),
\]
so that $\pi_c^{\Gamma}(M)\in\mbox{\rm Mat}_{c\times c}^{\Gamma}$.
Moreover, we clearly have that $\pi_c^{\Gamma}(M)=M$ if and only if
$M\in\mbox{\rm Mat}^{\Gamma}_{c\times c}$.
We then conclude that
$\pi_c^{\Gamma}$ is idempotent, with range $\mbox{\rm
  Mat}^{\Gamma}_{c\times c}$.  Finally, (\ref{morphism}) follows from 
the definition of $\pi_c^{\Gamma}$ and the fact that $A$ commutes with
$G(\gamma)$, for all $\gamma\in\Gamma$.
\hfill
\qed

\begin{prop}
Let ${\cal B}$ 
be a $p$-parameter unfolding of $B\in\mbox{\rm Mat}_{c\times c}^{\Gamma}$
such that (\ref{suff_cond_versal}) holds (which implies by
Proposition \ref{suff_cond_versal_prop} that ${\cal B}$ is a
versal unfolding of $B$).  
Let $\pi_c^{\Gamma}$ be the projection of $\mbox{\rm Mat}_{c\times c}$ onto
$\mbox{\rm Mat}_{c\times c}^{\Gamma}$ such as in Lemma
\ref{homlem}.  
Then ${\cal B}^{\Gamma}\equiv \pi_c^{\Gamma}({\cal B})$ is a $\Gamma$-versal
unfolding of $B$.
\label{rest_from_full}
\end{prop}
\proof 
We first show that $\pi_c^{\Gamma}(T\Sigma_B)=T\Sigma^{\Gamma}_B$.
Let $x=[B,y]\in T\Sigma^{\Gamma}_B$, where $y\in\mbox{\rm
 Mat}^{\Gamma}_{c\times c}$.
Since obviously $x\in\mbox{\rm Mat}^{\Gamma}_{c\times c}$,
then
$x=\pi_c^{\Gamma}(x)=\pi_c^{\Gamma}([B,y])\in\pi_c^{\Gamma}(T\Sigma_B)$,
since $[B,y]$ is also an element of $T\Sigma_B$.  Thus
$T\Sigma^{\Gamma}_B\subset \pi_c^{\Gamma}(T\Sigma_B)$.
Now let $x\in\pi_c^{\Gamma}(T\Sigma_B)$, i.e.
$x=\pi_c^{\Gamma}(By-yB)$ for some $y\in\mbox{\rm Mat}_{c\times c}$.
Then, from Lemma \ref{homlem}, we have
\[
x=\pi_c^{\Gamma}(By-yB)=\pi_c^{\Gamma}(By)-\pi_c^{\Gamma}(yB)=
[B,\pi_c^{\Gamma}(y)]\in T\Sigma^{\Gamma}_B.
\]
Therefore, we have shown that 
$\pi_c^{\Gamma}(T\Sigma_B)=T\Sigma^{\Gamma}_B$.

Now, we will show that $\pi_c^{\Gamma}(D_{\alpha}{\cal
  B}(\alpha_0)\cdot {\mathbb C}^p)=D_{\alpha}{\cal
  B}^{\Gamma}(\alpha_0)\cdot {\mathbb C}^p$.  We have that \mbox{\rm $x\in
D_{\alpha}{\cal B}^{\Gamma}(\alpha_0)\cdot {\mathbb C}^p$}, if and only
  if there
exists $v\in {\mathbb C}^p$ such that
\[
\begin{array}{lll}
x&=&{\displaystyle\left.\frac{d}{d\varepsilon}\left({\cal
      B}^{\Gamma}(\alpha_0+\varepsilon
      v)\right)\right|_{\varepsilon=0}=
\left.\frac{d}{d\varepsilon}\left(\pi_c^{\Gamma}({\cal
      B}(\alpha_0+\varepsilon
      v)\right)\right|_{\varepsilon=0}}\\[0.15in]
&=&{\displaystyle
\pi_c^{\Gamma}\left(\left.\frac{d}{d\varepsilon}\left({\cal
      B}(\alpha_0+\varepsilon
      v)\right)\right|_{\varepsilon=0}\right)\in
\pi_c^{\Gamma}(D_{\alpha}{\cal B}(\alpha_0)\cdot {\mathbb C}^p).}
\end{array}
\]
Therefore, $\pi_c^{\Gamma}(D_{\alpha}{\cal
  B}(\alpha_0)\cdot {\mathbb C}^p)=D_{\alpha}{\cal
  B}^{\Gamma}(\alpha_0)\cdot {\mathbb C}^p$.

Finally, we show that 
\[
\mbox{\rm Mat}^{\Gamma}_{c\times c}=T\Sigma^{\Gamma}_B+D_{\alpha}{\cal
  B}^{\Gamma}(\alpha_0)\cdot {\mathbb C}^p.
\]
Obviously, from the preceding results, we have that
\[
T\Sigma^{\Gamma}_B+D_{\alpha}{\cal
  B}^{\Gamma}(\alpha_0)\cdot {\mathbb C}^p\subset \mbox{\rm Mat}^{\Gamma}_{c\times c}.
\]
Let $x\in\mbox{\rm Mat}^{\Gamma}_{c\times c}$, then from 
(\ref{suff_cond_versal}) we have that $x=x_a+x_b$, where
$x_a\in T\Sigma_B$ and $x_b\in D_{\alpha}{\cal B}(\alpha_0)\cdot
{\mathbb C}^p$.  Therefore
\[
x=\pi_c^{\Gamma}(x)=\pi_c^{\Gamma}(x_a)+\pi_c^{\Gamma}(x_b)\in
\pi_c^{\Gamma}(T\Sigma_B)+\pi_c^{\Gamma}(D_{\alpha}{\cal
  B}(\alpha_0)\cdot {\mathbb C}^p)=
T\Sigma^{\Gamma}_B+D_{\alpha}{\cal B}^{\Gamma}(\alpha_0)\cdot {\mathbb
  C}^p.
\]
Thus, we have shown that $\mbox{\rm Mat}^{\Gamma}_{c\times c}=T\Sigma^{\Gamma}_B+D_{\alpha}{\cal
  B}^{\Gamma}(\alpha_0)\cdot {\mathbb C}^p$, and by Lemma
\ref{suff_cond_versal_prop_gamma}, we conclude that ${\cal
  B}^{\Gamma}$ is a $\Gamma$-versal unfolding of $B$.
\hfill
\qed

\begin{rmk}\label{rmk:mini}
 Let $B\in \mbox{\rm Mat}^{\Gamma}_{c\times c}$, be given, and let
${\cal B}(\alpha)$ be a $p$-parameter unfolding of $B$ which is mini-versal.
Then (\ref{suff_cond_versal}) holds as a direct sum, and by the previous
Proposition, ${\cal B}^{\Gamma}(\alpha)=\pi^{\Gamma}_c({\cal B}(\alpha))$ is
a $\Gamma$-versal unfolding of $B$.  It should be noted that in general,
${\cal B}^{\Gamma}(\alpha)$ need not be $\Gamma$-mini-versal, even though
${\cal B}$ is mini-versal (in fact it is quite easy to constrct examples
to illustrate this fact). However, it is always possible to extract a 
$\Gamma$-mini-versal unfolding of $B$ from ${\cal B}^{\Gamma}$, using a 
straighforward procedure as we now illustrate.
\end{rmk} 
Let ${\cal W}$ be any direct sum complement of $T_B\Sigma$ in $\mbox{\rm Mat}_{c\times c}$,
i.e.
\begin{equation}
\mbox{\rm Mat}_{c\times c}=T_B\Sigma\oplus {\cal W}.
\label{TBW}
\end{equation}
For example, if ${\cal W}$ is chosen to be the space which is spanned by the
matrices $\Omega_1,\ldots,\Omega_{p}$ described in section 6 of \cite{BL02},
then the decomposition of any $c\times c$ matrix using (\ref{TBW}) is
particularly simple to compute.
Now, for each $i=1,\ldots,p$, we define
\[
\hat{B}^{\Gamma}_i\equiv 
\frac{d}{d\varepsilon}\left.\left({\cal B}^{\Gamma}(\alpha_0+\varepsilon\,e_i)
\right)\right|_{\varepsilon=0},
\]
where $e_i$ is the $p$-dimensional vector whose $i^{th}$ component is 1 and
all other components are zero.
Note that
\[
D_{\alpha}{\cal B}^{\Gamma}(\alpha_0)\cdot {\mathbb C}^p=
{\mathbb C}\cdot\{\,\hat{B}^{\Gamma}_1,\ldots,\hat{B}^{\Gamma}_p\,\}.
\]
Now, decompose these matrices using (\ref{TBW}), i.e.
write
\begin{equation}\label{mattheta}
\hat{B}^{\Gamma}_i=[B,y_i]+\sum_{j=1}^p\,\theta_{ij}\Omega_j,\,\,\,\,i=1,\ldots,p,
\end{equation}
where $y_i\in \mbox{\rm Mat}_{c\times c}$, and the $\theta_{ij}$ are scalars.
Let $\Theta=(\theta_{ij})$ be the $p\times p$ matrix coming from the
decomposition~(\ref{mattheta}). If we then let rows $i_1,\ldots i_k$ of the 
matrix $\Theta$ determine a maximal set of linearly independent vectors in 
the rowspace of the matrix $\Theta$, it follows from simple linear 
algebra that
\[
\mbox{\rm Mat}^{\Gamma}_{c\times c}=T_B\Sigma^{\Gamma}\oplus {\mathbb C}\cdot
\{\,\hat{B}^{\Gamma}_{i_1},\ldots,\hat{B}^{\Gamma}_{i_k}\,\},
\]
and so the following $k$-parameter $\Gamma$-unfolding of $B$ is $\Gamma$-mini-versal:
\[
\hat{{\cal B}}^{\Gamma}(\beta_1,\ldots,\beta_k)\equiv
{\cal B}^{\Gamma}\left(\alpha_0+\sum_{\ell=1}^k\,\beta_{\ell}\,e_{i_{\ell}}\right).
\]

Recall that the RFDE (\ref{fde}) is defined on ${\mathbb C}^n$, and that
$\Gamma$ is a subgroup of $GL(n,{\mathbb C})$.  It will be useful to
define a projection analogous to $\pi_c^{\Gamma}$ on $\mbox{\rm
  Mat}_{n\times n}$.
\begin{Def}
Let $A\in\mbox{\rm Mat}_{n\times n}$.  We define
\[
\pi_n^{\Gamma}(A)=\int_{\Gamma}\,\gamma\,A\,\gamma^{-1}\,d\gamma.
\]
Obviously, it follows that for all
$A\in\mbox{\rm Mat}_{n\times n}$, we have
$\sigma\, \pi^{\Gamma}_n(A)=
\pi^{\Gamma}_n(A)\,\sigma$ for all $\sigma\in\Gamma$.
\label{nprojdef}
\end{Def}

\begin{prop}
Suppose $C\in\mbox{\rm Mat}_{c\times c}$ is such that
\begin{equation}\label{Ceq}
C=\Psi(0)A\Phi(\tau^*),
\end{equation}
where $A$ is some $n\times n$ matrix, and $\tau^*$ is some fixed
number in $[-\tau,0]$.  Then
\[
\pi^{\Gamma}_c(C)=\Psi(0)\pi^{\Gamma}_n(A)\Phi(\tau^*).
\]
\label{Aprop}
\end{prop}
\proof
This is just a simple computation, using (\ref{c_act}):
\[
\pi^{\Gamma}_c(C)=\int_{\Gamma}\,G(\gamma)\Psi(0)A\Phi(\tau^*)G(\gamma^{-1})\,d\gamma
= \Psi(0)\left(\int_{\Gamma}\,\gamma\,A\,\gamma^{-1}\,d\gamma\right)\Phi(\tau^*).
\]
\hfill
\qed

\Section{Construction of a {\boldmath $\Gamma$}-equivariant {\boldmath $\Lambda$}-versal unfolding}

Consider now a $C^{\infty}$ parametrized family of 
linear RFDEs
of the form
\begin{equation}
\dot{z}(t)={\cal L}(\alpha)(z_t),
\label{fdep}
\end{equation}
where $\alpha\in {\mathbb C}^p$, and ${\cal L}(\alpha_0)={\cal L}_0$ is as in 
(\ref{fde}), for some $\alpha_0\in {\mathbb C}^p$.  
In the sequel, we will assume that a translation has been 
performed in the parameter space ${\mathbb C}^p$ such that
$\alpha_0=0$.
It was shown in \cite{BL02} that for a given set $\Lambda$ of 
solutions to (\ref{char_eq}),
a parameter-dependent reduction of (\ref{fdep}) to the $c$-dimensional
generalized
eigenspace corresponding to $\Lambda$ is given by
\begin{equation}
\dot{x}=\tilde{{\cal B}}(\alpha)x, 
\end{equation}
where 
\begin{equation}
\tilde{{\cal B}}(\alpha)=B+\Psi(0)[{\cal L}(\alpha)-{\cal L}_0](\Phi+
h(\alpha)),
\label{Bunfdef}
\end{equation}
where $h$ is a smooth $n\times c$ matrix-valued
function such that $h(0)=0$, and $B$ is as in
(\ref{red_ode}). 
\begin{rmk}
The function $h$ is the unique solution (via the implicit function
theorem) of a certain nonlinear equation $N(\alpha,h)=0$ (see equation
(4.5) in Proposition 4.1 of \cite{BL02}).  It is a straightforward
computation that for all $h$ and all $\gamma\in\Gamma$, we have
\[
N(\alpha,\gamma\,h\,G(\gamma^{-1}))=\gamma\,N(\alpha,h)\,G(\gamma^{-1})
\]
so that if $h(\alpha)$ is such that $N(\alpha,h(\alpha))=0$, then by
uniqueness of the solution, we must have
$\gamma\,h(\alpha)\,G(\gamma^{-1})=h(\alpha)$, i.e.
\[
\gamma\,h(\alpha)=h(\alpha)G(\gamma^{-1})
\]
 for all $\alpha$ near
$0$ and for all $\gamma\in\Gamma$.  Thus, if (\ref{fdep}) is
$\Gamma$-equivariant for all $\alpha$, 
then it follows from (\ref{Bunfdef}) that
$\tilde{{\cal B}}(\alpha)\in\mbox{\rm Mat}_{c\times c}^{\Gamma}$ for all
$\alpha$ near 0 in ${\mathbb C}^p$, i.e. $\tilde{{\cal B}}$ is a
$\Gamma$-unfolding of $B$.
\label{rmkequivh}
\end{rmk}
 
We recall the following definition from \cite{BL02}.
\begin{Def}
The parametrized family of RFDEs~(\ref{fdep}) is said to be a 
{\em $\Lambda$-versal unfolding} (respectively a {\em $\Lambda$-mini-versal
unfolding}) for
the RFDE~(\ref{fde}) if the matrix $\tilde{{\cal B}}(\alpha)$ defined by~(\ref{Bunfdef}) is 
a versal unfolding (respectively a mini-versal unfolding) for $B$.
\end{Def}
In the case where (\ref{fde}) is $\Gamma$-equivariant, we have
\begin{Def}
The parametrized family (\ref{fdep}) is said to be a
{\em $\Gamma$-equivariant $\Lambda$-versal unfolding} (respectively a
{\em $\Gamma$-equivariant $\Lambda$-mini-versal unfolding}) for
(\ref{fde})
if (\ref{fdep}) is $\Gamma$-equivariant for all $\alpha$, and the
matrix
$\tilde{{\cal B}}(\alpha)$ defined by (\ref{Bunfdef}) is a $\Gamma$-versal
unfolding
(respectively a $\Gamma$-mini-versal unfolding) for $B$.
\end{Def}

In \cite{BL02}, we showed that given $B\in\mbox{\rm Mat}_{c\times c}$,
there exists points $\tau_0,\tau_1,\ldots,\tau_{c-q}$ in $[-\tau,0]$,
and $n\times n$ matrices $A^m_j$, $j=0,\ldots,c-q$,
$m=0,\ldots,\delta$
(where $\delta$ is the codimension of $T\Sigma_B$ in $\mbox{\rm
  Mat}_{c\times c}$) such that the following $\delta$-parameter
unfolding
of $B$
\begin{equation}
{\cal
  B}(\alpha)=B+\sum_{m=1}^{\delta}\,\alpha_m\,\left(\sum_{j=0}^{c-q}\,\Psi(0)A^m_j\Phi(\tau_j)\right)
\label{paper1unfdef}
\end{equation}
satisfies (\ref{suff_cond_versal}) as a direct sum (i.e. ${\cal B}$ is
a
mini-versal unfolding of $B$).  From this, it follows 
(see Theorem 7.4 of \cite{BL02})      
that if ${\cal L}(\alpha)$ is 
the $\delta$-parameter family of bounded linear
operators from $C([-\tau,0],{\mathbb C}^n)$ into ${\mathbb C}^n$ defined in the
following way:
\[
{\cal L}(\alpha)\,z={\cal
  L}_0\,z+\sum_{m=1}^{\delta}\,\alpha_m\,\left(\sum_{j=0}^{c-q}\,A^m_j\,z(\tau_j)\right),
\]
where the $\alpha_m$ are complex parameters and ${\cal L}_0$ is as in 
(\ref{fde}), then (\ref{fdep}) is a $\Lambda$-versal unfolding of (\ref{fde}).
 
We now state and prove the main result of this paper.
\begin{thm}
Suppose $B\in\mbox{\rm Mat}^{\Gamma}_{c\times c}$, and
consider the matrices $A^m_j$ which appear in
the versal unfolding (\ref{paper1unfdef}) of $B$.  
Then the family ${\cal L}(\alpha)$ of $\Gamma$-equivariant bounded linear operators
defined by
\begin{equation}
{\cal L}(\alpha)\,z={\cal L}_0\,z+\sum_{m=1}^{\delta}\,\alpha_m\,\left(
\sum_{j=0}^{c-q}\,\pi^{\Gamma}_n(A^m_j)\,z(\tau_j)\right),
\label{unf_equiv}
\end{equation}
(where $\pi^{\Gamma}_n$ is as in Definition \ref{nprojdef}) is such that
(\ref{fdep}) is a $\Gamma$-equivariant $\Lambda$-versal unfolding of
(\ref{fde}).
\label{main_thm}
\end{thm}
\proof
From Propositions \ref{rest_from_full} and \ref{Aprop}, we conclude
that the versal unfolding (\ref{paper1unfdef}) of $B$ is such that
\[
\pi^{\Gamma}_c({\cal B}(\alpha))=B+
\sum_{m=1}^{\delta}\,\alpha_m\,\left(\sum_{j=0}^{c-q}\,\Psi(0)\pi^{\Gamma}_n(A^m_j)\Phi(\tau_j)\right)
\]
is a $\Gamma$-versal unfolding of $B$.  Thus, if we define 
${\cal L}(\alpha)$ as in (\ref{unf_equiv}), then $\tilde{{\cal
    B}}(\alpha)$
(as in (\ref{Bunfdef})) is 
such that
\[
\tilde{{\cal B}}(\alpha)=\pi^{\Gamma}_c({\cal B}(\alpha))+
\Psi(0)[{\cal L}(\alpha)-{\cal L}_0](h(\alpha)).
\]
We conclude that $\tilde{{\cal B}}(\alpha)$ is
also a $\Gamma$-versal unfolding of $B$,
since $h(0)=0$ implies that $D_{\alpha}\tilde{{\cal B}}(0)=
D_{\alpha}\pi^{\Gamma}_c({\cal B}(\alpha))|_{\alpha=0}$.
\hfill
\qed

\Section{An example with {\boldmath $\Gamma=\D_3$}}
In this section we compute the versal unfolding of a $\D_3$-equivariant
system of linear delay-differential equations near points of double Hopf
bifurcations.
Consider the group $\D_3$ generated by $\{\kappa,\gamma\}$
with $\kappa^2=\gamma^3=1$ and the representation $\rho:\D_3\to GL(3,{\mathbb C})$ 
given by
\begin{equation}\label{act3}
\rho(\kappa)=\left(\begin{array}{ccc}
1 & 0 & 0 \\
0 & 0 & 1 \\
0 & 1 & 0 
\end{array}\right),\quad 
\rho(\gamma)=\left(\begin{array}{ccc}
0 & 1 & 0 \\
0 & 0 & 1 \\
1 & 0 & 0 
\end{array}\right).
\end{equation}

Wu {\em et al.}~\cite{WFH99} and Ncube {\em et al.}~\cite{NCW} consider the 
$\D_3$-equivariant linear delay-differential equation
\begin{equation}\label{d3eqs}
\begin{array}{rcl}
\dot u_{j}(t)&=&-u_{j}(t)+\alpha u_{j}(t-\tau_{s})
+\beta\left[u_{j-1}(t-\tau_n)+u_{j+1}(t-\tau_{n})\right],\quad j=1,2,3
\end{array}
\end{equation}
where the indices are taken $\mod 3$. This linear system appears as the 
linearization at an equilibrium of a coupled system of three multiple-delayed 
identical neurons. In those papers, the focus is on finding periodic
solutions via the Hopf bifurcation theorem.

The spectrum of the linear operator associated to~(\ref{d3eqs}) is obtained by solving the 
characteristic equation
\[
S(\lambda)=\Delta_1(\lambda)\Delta_2(\lambda)^{2}=0,
\]
where 
\[
\begin{array}{rcl}
\Delta_1(\lambda)&=&\lambda+1-\alpha e^{-\lambda\tau_s}
-2\beta e^{-\lambda\tau_n},\\
\Delta_2(\lambda)&=&\lambda+1-\alpha e^{-\lambda\tau_s}
+\beta e^{-\lambda\tau_n}.
\end{array} 
\]
We now sketch an argument showing that points of double Hopf
bifurcation exist for both factors $\Delta_1$ and $\Delta_2$ of 
the characteristic equation. 

\begin{figure}
\centerline{%
(a) \psfig{file=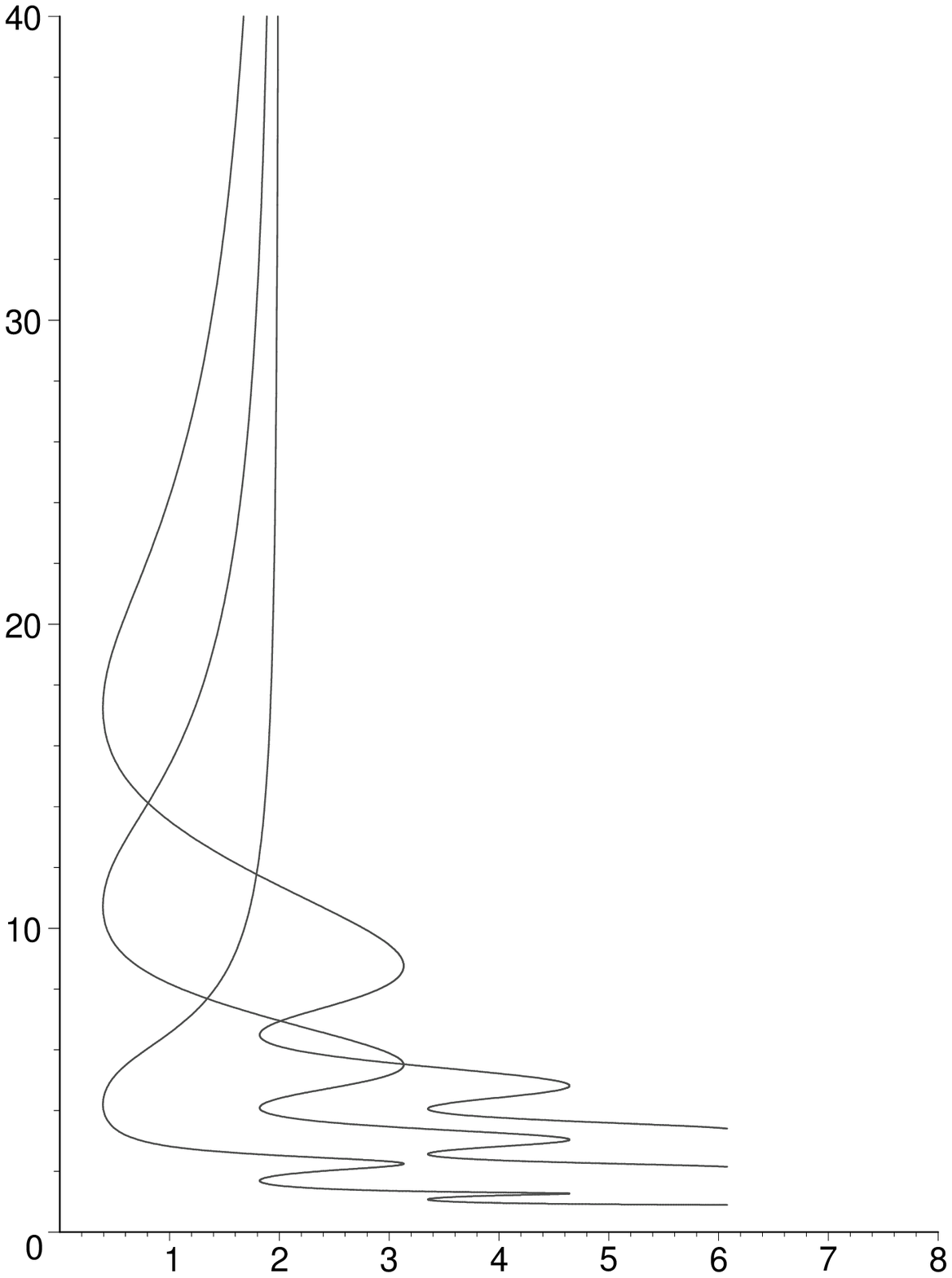,width=2.0in,height=2.0in}
\qquad
(b) \psfig{file=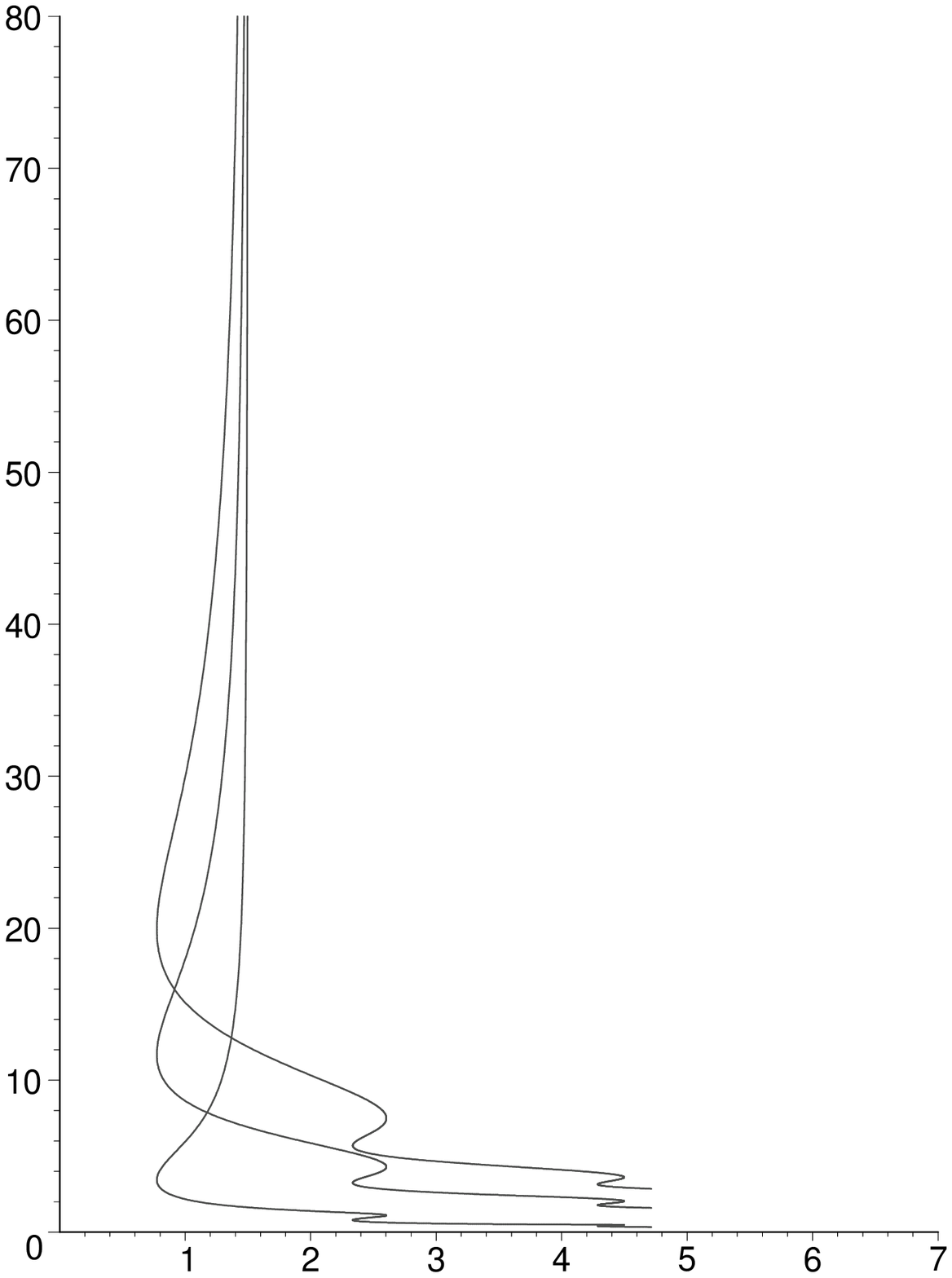,width=2.0in,height=2.0in}}
\caption{Curves of imaginary eigenvalues in $(\alpha,\tau_s)$-space
for (a) $\Delta_1=0$ where $\beta=-0.5$ and $\tau_n=4$, and
(b) $\Delta_2=0$ where $\beta=0.5$ and $\tau_n=3$.}
\label{2hopf}
\end{figure}

Suppose that we are looking for a point of double Hopf
bifurcation solution to $\Delta_1=0$, then the eigenvalues 
$i\omega$ must satisfy
\[
\begin{array}{rcl}
\omega+\alpha \sin(\omega\tau_s)&=&-2\beta\sin (\omega\tau_n)\\
1-\alpha\cos(\omega\tau_s)&=&2\beta\cos(\omega\tau_n).
\end{array}
\]
It is easy to solve these equations for $\alpha$ and $\tau_s$
in terms of the other parameters, we obtain the following two
equations 
\[
\begin{array}{rcl}
\alpha&=&\pm\sqrt{(1-2\beta\cos(\omega\tau_n))^{2}
+(\omega+2\beta\sin(\omega\tau_n))^{2}}\\
\tau_s&=&\dfrac{1}{\omega}\arctan\left(\dfrac{-\omega
-2\beta\sin(\omega\tau_n)}{1-2\beta\cos(\omega\tau_n)}\right).
\end{array}
\]
Then $\alpha$ and $\tau_s$ with $\beta$ and $\tau_n$ fixed represent 
curves parametrized by $\omega$ in $(\alpha,\tau_s)$-space. Note that 
there are several branches of curves depending on the branch of the
arctan which is chosen. Figure~\ref{2hopf}(a) shows some of these 
curves for $\beta=-0.5$ and $\tau_n=4$. Points of double 
Hopf bifurcation lie at the intersection of these curves.
Points of double Hopf bifurcation for $\Delta_2=0$ are found in
a similar way by isolating $\alpha$ and $\tau_s$ as functions
of $\beta$, $\tau_n$ and $\omega$. Figure~\ref{2hopf}(b) shows 
several curves of imaginary eigenvalues plotted for $\beta=0.5$ 
and $\tau_n=3$.

From the graphics in Figure~\ref{2hopf}, it is reasonable to suppose 
that the characteristic equation has a point of double Hopf bifurcation at  
$(\alpha^{*},\beta^{*},\tau_s^{*},\tau_n^{*})$ with $\tau_s^{*}$, $\tau_n^{*}$
nonzero, $\tau_s^{*}\neq \tau_n^{*}$ and frequencies $\omega_1\neq \omega_2$. 
Of course, the point $(\alpha^{*},\beta^{*},\tau_s^{*},\tau_n^{*})$ as 
well as the frequencies are different whether $\Delta_1=0$ or $\Delta_2=0$,
however, to keep the notation to a minimum we use this unique label for 
both cases.

We let $\Lambda=\{\pm i\omega_1,\pm i\omega_2\}$ with $\omega_1\neq 0$ and
$\omega_2\neq 0$ be the set of eigenvalues on the imaginary axis at 
a double Hopf bifurcation point.
If the double Hopf point occurs for the $\Delta_1$ factor, then
the eigenvalues of $\Lambda$ are simple and the restriction 
of the linear operator to the center eigenspace is given by
\[
B=\left(\begin{array}{cccc}
\omega_1 i & 0 & 0 & 0 \\
 0 & -\omega_1 i & 0 & 0 \\
 0 & 0           & \omega_2 i & 0 \\
 0 & 0           & 0          & -\omega_2 i
\end{array}\right).
\]
If the double Hopf point occurs for the $\Delta_2$ factor, 
then the eigenvalues are double and the linear flow on the
center eigenspace is given by the following matrix
\[
B=\left(\begin{array}{cccccccc}
\omega_1 i & 0 & 0 & 0 & 0 & 0 & 0 & 0\\
 0 & \omega_1 i & 0 & 0 & 0 & 0 & 0 & 0\\
 0 & 0 & -\omega_1 i & 0 & 0 & 0 & 0 & 0\\
 0 & 0 & 0 & -\omega_1 i & 0 & 0 & 0 & 0 \\
 0 & 0 & 0 & 0 & \omega_2 i & 0 & 0 & 0 \\
 0 & 0 & 0 & 0 & 0 & \omega_2 i & 0 & 0 \\
 0 & 0 & 0 & 0 & 0 & 0 & -\omega_2 i & 0 \\
 0 &  0 & 0 & 0 & 0 & 0 & 0 & -\omega_2 i
\end{array}\right).
\]
We now compute the $\D_3$-equivariant $\Lambda$-mini-versal 
unfolding of the linear bounded operator in both cases.
Recall that our goal is to find matrices $A_{j}^{m}$ from which we
can write~(\ref{unf_equiv}). The computations rely on 
the results in~\cite{BL02}, mainly Section 6 and Lemma 7.1.

\subsection{Simple imaginary eigenvalues.}
The center eigenspace $P$ has basis
\[
\begin{array}{rcl}
\Phi(\theta)&=&(\phi(\theta),\phi_2(\theta),\phi_3(\theta),\phi_4(\theta))\\
&=&(u e^{i\omega_1 \theta},u e^{-i\omega_1 \theta},u e^{i\omega_2 \theta},u e^{-i\omega_2 \theta})
\end{array}
\]
where
\[
u=\left(\begin{array}{c}
1\\ 1 \\ 1 
\end{array}\right).
\]
Its dual, $P^{*}$, has basis $\Psi(s)=\left(\psi_1(s),\psi_2(s),\psi_3(s),
\psi_4(s)\right)^{t}$ where $\psi_2=\ov{\psi_1}$, $\psi_4=\ov{\psi_3}$ and 
$\psi_{j}(s)=(\psi_{j1}(s),\psi_{j2}(s),\psi_{j3}(s))$, $j=1,2,3,4$.

\paragraph{Basis of the unfolding space.}
Let ${\cal W}$ be the unfolding space spanned by matrices $\Omega_1,\ldots,
\Omega_p$ defined in Section 6 of~\cite{BL02}. The dimension of ${\cal W}$ 
for matrix $B$ is given by
\begin{equation}\label{dimW}
\sum_{j=1}^{r}\sum_{\ell=1}^{k_{j}}(2\ell-1)n_{j,\ell}
\end{equation}
where $r=4$ is the number of distinct eigenvalues, $k_{j}=1$
is the number of Jordan blocks for each eigenvalue and $n_{j,\ell}=1$
is the size of each Jordan block, thus the unfolding space is 
$4$-dimensional. The space ${\cal E}({\cal W})$ lying in the range of
$\Psi(0)$ is linearly isomorphic to ${\cal W}$ and is also an unfolding 
space, see Section 6 of~\cite{BL02}. A basis of ${\cal E}({\cal W})$ is 
given by $\Psi(0)R^{j}$ for $j=1,2,3,4$ where  $R^1$, $R^2$, $R^3$, $R^4$ 
are $3\times 4$ matrices defined following the procedure established 
in~\cite{BL02}, Section 6.3. In~\cite{BL02} we introduce
projections $\Pi_{j}:\C^{3}\to \C$ ($j=1,2,3,4)$ with
\[
\Pi_{j}(v)=\psi_{j1}(0) v^1+\psi_{j2}(0) v^2+\psi_{j3}(0)v^3.
\]
Let $v_{j}=(v_{j}^1,v_{j}^2,v_{j}^3)^{t}\in\C^{3}$ ($j=1,2,3,4$) be 
such that $\Pi_{j}(v_{j})=1$ and set $R^1=(v_1,0,0,0)$, $R^2=(0,v_2,0,0)$,
$R^3=(0,0,v_3,0)$, and $R^{4}=(0,0,0,v_4)$. Note that since 
$\psi_2=\ov{\psi_1}$ and $\psi_4=\ov{\psi_3}$, we choose 
$v_2=\ov{v_1}$ and $v_4=\ov{v_3}$.

Since $\mbox{rank}\,\Phi(0)=1$, by Theorem 7.2 of~\cite{BL02} we need 
$4$ distinct delays $\tau_0,\tau_1,\tau_2,\tau_3$ to find $3\times 3$
matrices $A_{j}^{i}$ solving
\begin{equation}\label{Reqs}
R^{i}=\sum_{j=0}^{3} A_{j}^{i}\Phi(\tau_{j}).
\end{equation}

\paragraph{Projection of the basis.}
By Theorem~\ref{main_thm}, the solutions to~(\ref{Reqs}) are used to 
write the equivariant $\Lambda$-unfolding by projecting the matrices 
$A_{j}^{m}$ with $\pi_{n}^{\Gamma}$. Multiplying~(\ref{Reqs}) by $\Psi(0)$ 
we obtain equation~(\ref{Ceq}) with $C=\Psi(0)R$. Then, projecting $C$ 
with $\pi_{c}^{\Gamma}$, we obtain
\begin{equation}\label{newproj}
\pi_{c}^{\Gamma}(C)=\dps\int_{\Gamma} G(\gamma)\Psi(0)R\,G(\gamma^{-1})d\gamma
=\dps\int_{\Gamma} \Psi(0)\cdot\rho(\gamma)R\, G(\gamma^{-1})d\gamma=\Psi(0)\dps\int_{\Gamma} \rho(\gamma)R\, G(\gamma^{-1})d\gamma.
\end{equation}
We project the matrices $R^{j}$ on the subspace of $\D_3$-equivariant 
$3\times 4$ matrices using the projection operator obtained 
from~(\ref{newproj}); that is, since the $G$ representation is trivial 
on $P$ and $\Gamma=\D_3$ is finite, then the projected matrices are
\[
\dps\ov{R}^j\equiv\dfrac{1}{6}\sum_{g\in \D_3} \rho(g) R^{j}.
\]
Note that the factor $\frac{1}{6}$ is not necessary to the 
projection of the matrices and could be removed without further
consequences.
Thus, we obtain new matrices
\[
\ov{R}^1=\left(\begin{array}{cccc}
\eta_1 & 0 & 0 & 0\\
\eta_1 & 0 & 0 & 0\\
\eta_1 & 0 & 0 & 0
\end{array}\right), \quad
\ov{R}^2=\left(\begin{array}{cccc}
0 & \eta_2 & 0 & 0 \\
0 & \eta_2 & 0 & 0 \\
0 & \eta_2 & 0 & 0 
\end{array}\right),
\]
\[
\ov{R}^3=\left(\begin{array}{cccc}
0 & 0 & \eta_3 & 0\\
0 & 0 & \eta_3 & 0\\
0 & 0 & \eta_3 & 0
\end{array}\right), \quad
\ov{R}^4=\left(\begin{array}{cccc}
0 & 0 & 0 & \eta_4\\
0 & 0 & 0 & \eta_4\\
0 & 0 & 0 & \eta_4 
\end{array}\right)
\]
where $\eta_j=\frac{1}{3}(v_{j}^{1}+v_{j}^{2}+v_{j}^{3})$ and of course
$\eta_2=\ov{\eta_1}$ and $\eta_4=\ov{\eta_3}$. Because the matrices
$\ov{R}^{j}$ are $\D_3$-equivariant matrices, we may take matrices 
$A_{j}^{m}$ $\D_3$-equivariant with respect to the action of $\D_3$ 
on $\C^{3}$ given by~(\ref{act3}); that is, 
\begin{equation}\label{D3mat3}
A_{j}^{m}=\left(\begin{array}{ccc}
a_j^m & b_j^m & b_j^m \\
b_j^m & a_j^m & b_j^m \\
b_j^m & b_j^m & a_j^m
\end{array}\right).
\end{equation}

\paragraph{$\D_3$-mini-versal unfolding}
Before solving equation~(\ref{Reqs}), we verify which of the four
matrices $\ov{R}^{j}$ ($j=1,2,3,4$) are necessary to obtain 
a mini-versal unfolding. We follow the steps described after 
Remark~\ref{rmk:mini}.
Since ${\cal E}({\cal W})$ is an unfolding space for $B$, we 
set $\hat{B}_{j}^{\D_3}=\Psi(0)\ov{R}^{j}$ for $j=1,2,3,4$. 
The space ${\cal W}$ is spanned by $4\times 4$ matrices
$\Omega_i$ ($i=1,2,3,4$) where $\Omega_i$ has a $1$ on the diagonal 
at $(i,i)$ and is zero elsewhere. Then, it is easy to compute that
\[
\begin{array}{rcl}
\hat{B}_1^{\D_3}&=&[B,y_1]+\eta_1(\psi_{11}(0)+\psi_{12}(0)+\psi_{13}(0))
\Omega_1\\[0.15in]
\hat{B}_2^{\D_3}&=&[B,y_2]+\eta_2(\ov{\psi_{11}(0)}+\ov{\psi_{12}(0)}+
\ov{\psi_{13}(0)})\Omega_2\\[0.15in]
\hat{B}_3^{\D_3}&=&[B,y_3]+\eta_3(\psi_{31}(0)+\psi_{32}(0)+\psi_{33}(0))\Omega_3\\[0.15in]
\hat{B}_4^{\D_3}&=&[B,y_4]+\eta_4(\ov{\psi_{31}(0)}+\ov{\psi_{32}(0)}+\ov{\psi_{33}(0)})\Omega_4
\end{array}
\]
for some matrices $y_1,y_2,y_3,y_4\in \mbox{Mat}_{4\times 4}$. 
Since the matrix $\Theta$ (defined using~(\ref{mattheta})) is diagonal, 
the four rows are needed to span the rowspace of $\Theta$. Hence,
\[
\mbox{Mat}_{4\times 4}^{\D_3}=T_{B}\Sigma^{\D_3}\oplus 
\C\cdot\{\hat{B}_1^{\D_3},\hat{B}_2^{\D_3},\hat{B}_3^{\D_3},
\hat{B}_4^{\D_3}\}
\]
and we need to consider all matrices $\ov{R}^{j}$ ($j=1,2,3,4$)
in order to compute the $\D_3$-equivariant $\Lambda$-mini-versal unfolding.

\paragraph{Computation of the matrices $\mathbf{A}_{j}^{m}$.}
We now solve for the matrices $A_{j}^{m}$ using equation~(\ref{Reqs}).
Since 
\[
A_{j}^{m}\Phi(\tau_{j})=(a_j^m+2b_j^m)\Phi(\tau_{j}),
\]
it is sufficient to solve the first line of~(\ref{Reqs}) to obtain
the matrices $A_j^m$. Consider the first line of equation~(\ref{Reqs}) 
for $\ov{R}^1$ 
\[
(\eta_1,0,0,0)=\sum_{j=0}^{3} (a_j^{1},b_{j}^{1},b_{j}^{1})\Phi(\tau_j).
\]
This equation can be rewritten as
\begin{equation}\label{ajieqs}
\left(\begin{array}{c}
\eta_1\\
0\\
0\\
0
\end{array}\right)=M
\left(
\begin{array}{c}
a_0^{1}+2 b_{0}^{1}\\ a_1^{1}+2 b_{1}^{1}\\ a_2^{1}+2 b_{2}^{1}\\ 
a_3^{1}+2 b_{3}^{1}
\end{array}
\right),
\end{equation}
where
\begin{equation}\label{matM}
M=\left(\begin{array}{cccc}
e^{i\omega_1 \tau_0} & e^{i\omega_1 \tau_1} & e^{i\omega_1 \tau_2} &
e^{i\omega_1 \tau_3} \\
e^{-i\omega_1 \tau_0} & e^{-i\omega_1 \tau_1} & e^{-i\omega_1 \tau_2} &
e^{-i\omega_1 \tau_3} \\
e^{i\omega_2 \tau_0} & e^{i\omega_2 \tau_1} & e^{i\omega_2 \tau_2} &
e^{i\omega_2 \tau_3} \\
e^{-i\omega_2 \tau_0} & e^{-i\omega_2 \tau_1} & e^{-i\omega_2 \tau_2} &
e^{-i\omega_2 \tau_3} 
\end{array}\right).
\end{equation}
Obviously, $\det M=0$ if any two of the delays are equal or if 
$\omega_1=\omega_2$. Now, since we take our four delays to be distinct
and it is assumed that $\omega_1\neq \omega_2$, then generically $M$ is
nonsingular. Thus, we obtain $L_1=\sum_{j=0}^{3} A_j^1 z(\tau_j)$ by 
solving~(\ref{ajieqs}). Now, the equation for $\ov{R}^{2}$ yields,
\[
\left(\begin{array}{c}
0\\
\ov{\eta}_1\\
0\\
0
\end{array}\right)=M
\left(
\begin{array}{c}
a_0^{2}+2 b_{0}^{2}\\ a_1^{2}+2 b_{1}^{2}\\ a_2^{2}+2 b_{2}^{2}\\ 
a_3^{2}+2 b_{3}^{2}
\end{array}
\right),
\]
which if we interchange rows $1$ with $2$ and $3$ with $4$ corresponds
to the complex conjugate of system~(\ref{ajieqs}). So, we can choose 
$a_j^2$ and $b_j^2$ so that $L_2=\ov{L}_1$.
In a similar manner we obtain $L_3$ and $L_4=\ov{L}_3$. 
 
Looking at equation~(\ref{ajieqs}), we see that for each $j=0,1,2,3$, 
either $a_j^1$ or $b_j^1$ can be chosen arbitrarily. We exploit this 
freedom to choose constants $a_{j}^{m}$ and $b_{j}^{m}$ which preserve 
the structure of equation~(\ref{d3eqs}). 

\paragraph{$\D_3$-equivariant $\Lambda$-unfolding.}
Take $\tau_0=0$, $\tau_1=\tau_s^{*}$, $\tau_2=\tau_n^{*}$, $b_{0}^{m}=0$,
$b_{1}^{m}=0$ and $a_{2}^{m}=0$ for $m=1,2,3,4$ (recall that we have assumed
$\tau_s^{*},\tau_n^{*}$ are nonzero and distinct so that matrix $M$ 
in~(\ref{matM}) is generically nonsingular with this substitution). Hence,
\[
\begin{array}{rcl}
L_1\,z&=&a_0^1 I z(0)+a_1^1 Iz(\tau_s^{*})+\left(\begin{array}{ccc} 
0 & b_2^1 & b_2^1 \\
b_2^1 & 0 & b_2^1 \\
b_2^1 & b_2^1 & 0
\end{array}\right) z(\tau_n^{*})+a_3^1 I z(\tau_3)+\left(\begin{array}{ccc} 
0 & b_3^1 & b_3^1 \\
b_3^1 & 0 & b_3^1 \\
b_3^1 & b_3^1 & 0
\end{array}\right) z(\tau_3)\\
L_3\,z&=&a_0^3 I z(0)+a_1^3 I z(\tau_s^{*})+\left(\begin{array}{ccc} 
0 & b_2^3 & b_2^3 \\
b_2^3 & 0 & b_2^3 \\
b_2^3 & b_2^3 & 0
\end{array}\right) z(\tau_n^{*})
+a_3^3 I z(\tau_3)+\left(\begin{array}{ccc} 
0 & b_3^3 & b_3^3 \\
b_3^3 & 0 & b_3^3 \\
b_3^3 & b_3^3 & 0
\end{array}\right) z(\tau_3),
\end{array}
\]
and the complex $\D_3$-equivariant $\Lambda$-versal unfolding is defined by
\[
\begin{array}{rcl}
{\cal L}(\alpha)\,z&=&{\cal L}_{0}\,z+\dps\sum_{m=1}^{4} \alpha_m L_m\\
&=&{\cal L}_{0}\,z+\dps\sum_{m=1}^{4} \alpha_{m} a_0^m I z(0)
+\dps\sum_{m=1}^{4} \alpha_{m} a_1^{m} I z(\tau_s^{*})+
\dps\sum_{m=1}^{4} \alpha_{m} \left(\begin{array}{ccc}
0 & b_2^m & b_2^m\\
b_2^m & 0 & b_2^m\\
b_2^m & b_2^m & 0  
\end{array}\right) z(\tau_n^{*})
\\
&+&\dps\sum_{m=1}^{4} \alpha_{m} a_3^m I z(\tau_3)+
\dps\sum_{m=1}^{4} \alpha_{m} \left(\begin{array}{ccc}
0 & b_3^m & b_3^m\\
b_3^m & 0 & b_3^m\\
b_3^m & b_3^m & 0  
\end{array}\right) z(\tau_3).
\end{array}
\]
By the form of equation~(\ref{ajieqs}), we can simplify the unfolding even
more by setting for instance $a_{3}^{m}=0$ for $m=1,2,3,4$. 
Thus a $\D_3$-equivariant $\Lambda$-versal unfolding of~(\ref{d3eqs})
near $(\alpha^{*},\beta^{*},\tau_s^{*},\tau_n^{*})$ is given by
\begin{equation}\label{unfDDE}
\begin{array}{rcl}
\dot u_{1}(t)&=&(-1+\epsilon_1)u_{1}(t)+(\alpha^{*}+\epsilon_2)u_{1}(t-\tau_s^{*})+(\beta^{*}+\epsilon_3)(u_{3}(t-\tau_n^{*})+u_{2}(t-\tau_n^{*}))\\
&+&\epsilon_4(u_{3}(t-\tau_3)+u_{2}(t-\tau_3))\\
\\
\dot u_{2}(t)&=&(-1+\epsilon_1)u_{2}(t)+(\alpha^{*}+\epsilon_2)u_{2}(t-\tau_s^{*})+(\beta^{*}+\epsilon_3)(u_{1}(t-\tau_n^{*})+u_{3}(t-\tau_n^{*}))\\
&+&\epsilon_4(u_{1}(t-\tau_3)+u_{3}(t-\tau_3))\\
\\
\dot u_{3}(t)&=&(-1+\epsilon_1)u_{3}(t)+(\alpha^{*}+\epsilon_2)u_{3}(t-\tau_s^{*})+(\beta^{*}+\epsilon_3)(u_{2}(t-\tau_n^{*})+u_{1}(t-\tau_n^{*}))\\
&+&\epsilon_4(u_{2}(t-\tau_3)+u_{1}(t-\tau_3)).
\end{array}
\end{equation}
where $\epsilon_1=\sum_{m=1}^{4} \alpha_{m} a_{0}^{m}$, 
$\epsilon_2=\sum_{m=1}^{4} \alpha_{m} a_{1}^{m}$,
$\epsilon_3=\sum_{m=1}^{4} \alpha_{m} b_{2}^{m}$ and
$\epsilon_4=\sum_{m=1}^{4} \alpha_{m} b_{3}^{m}$.

Now, the parameters $\epsilon_{i}$ are linearly independent
since
\begin{equation}\label{indep}
\begin{array}{rcl}
\left(\begin{array}{c}
\epsilon_1 \\ \epsilon_2\\ \epsilon_3\\ \epsilon_4
\end{array}\right)&=&
\left(\begin{array}{cccc}
a_0^1 & \ov{a}_0^1 & a_0^3 & \ov{a}_0^3\\
a_1^1 & \ov{a}_1^1 & a_1^3 & \ov{a}_1^3\\
b_2^1 & \ov{b}_2^1 & b_2^3 & \ov{b}_2^3\\
b_3^1 & \ov{b}_3^1 & b_3^3 & \ov{b}_3^3
\end{array}\right)
\left(\begin{array}{c}
\alpha_1 \\ \alpha_2 \\ \alpha_3 \\ \alpha_4
\end{array}\right)
=M^{-1} \,\mbox{diag}(\eta_1,\ov{\eta_1},\eta_3,\ov{\eta_3})
\left(\begin{array}{c}
\alpha_1 \\ \alpha_2 \\ \alpha_3 \\ \alpha_4
\end{array}\right).
\end{array}
\end{equation}
We now compute the real $\D_3$-equivariant $\Lambda$-versal unfolding.
From~\cite{BL02} we know that the real unfolding is given by
\[
{\cal L}(\alpha)={\cal L}_{0}+\alpha_{1} \re(L_1)+\alpha_2\im(L_1)+
\alpha_3 \re(L_3)+\alpha_4 \im(L_3)
\]
so the delay-differential equation is identical to~(\ref{unfDDE})
with  
\[
\begin{array}{rcl}
\epsilon_1&=&\alpha_1 \re(a_0^1)+\alpha_2 \im(a_0^1)+\alpha_3 \re(a_0^3)
+\alpha_4 \im(a_0^3),\\
\epsilon_2&=&\alpha_1 \re(a_1^1)+\alpha_2 \im(a_1^1)+\alpha_3 \re(a_1^3)
+\alpha_4 \im(a_1^3),\\
\epsilon_3&=&\alpha_1 \re(b_2^1)+\alpha_2 \im(b_2^1)+\alpha_3 \re(b_2^3)
+\alpha_4 \im(b_2^3),\\
\epsilon_4&=&\alpha_1 \re(b_3^1)+\alpha_2 \im(b_3^1)+\alpha_3 \re(b_3^3)
+\alpha_4 \im(b_3^3).
\end{array}
\]
Moreover, the new parameters $\epsilon_i$, $i=1,2,3,4$ are also linearly
independent since 
\[
\begin{array}{rcl}
\left(\begin{array}{c}
\epsilon_1 \\ \epsilon_2\\ \epsilon_3\\ \epsilon_4
\end{array}\right)&=&
\left(\begin{array}{cccc}
\re(a_0^1) & \im(a_0^1) & \re(a_0^3) & \im(a_0^3)\\
\re(a_1^1) & \im(a_1^1) & \re(a_1^3) & \im(a_1^3)\\
\re(b_2^1) & \im(b_2^1) & \re(b_2^3) & \im(b_2^3)\\
\re(b_3^1) & \im(b_3^1) & \re(b_3^3) & \im(b_3^3)
\end{array}\right)
\left(\begin{array}{c}
\alpha_1 \\ \alpha_2 \\ \alpha_3 \\ \alpha_4
\end{array}\right)
\end{array}
\]
where the matrix is nonsingular since the four columns of the matrix span 
the four dimensional real subspace of the four dimensional complex space 
spanned by the columns of the matrix in~(\ref{indep}).

\subsection{Double imaginary eigenvalues.}
We now consider the case where the eigenvalues of $\Lambda$ are double.
The space $P$ has basis
\[
\begin{array}{rcl}
\Phi(\theta)&=&(\phi_1(\theta),\phi_2(\theta),\phi_3(\theta),\phi_4(\theta),
\phi_5(\theta),\phi_6(\theta),\phi_7(\theta),\phi_8(\theta))\\
&=&(v e^{i\omega_1 \theta},\ov{v} e^{i\omega_1 \theta},\ov{v} e^{-i\omega_1 
\theta},v e^{-i\omega_1 \theta},v e^{i\omega_2 \theta},\ov{v} e^{i\omega_2 
\theta},\ov{v} e^{-i\omega_2 \theta},v e^{-i\omega_2 \theta}),
\end{array}
\]
where
\[
v=\left(\begin{array}{ccc}
1\\ \omega \\ \ov{\omega}
\end{array}\right)
\]
and $\omega=e^{2i\pi/3}$.
We now compute the basis $\Psi$ of the adjoint $P^{*}$. 
The action of $\D_3$ on $\C^{8}$ is generated by matrices
\begin{equation}\label{act8}
\begin{array}{rcl}
G(\gamma)&=&\mbox{diag}\,(e^{2\pi i/3},e^{-2\pi i/3},e^{-2\pi i/3},
e^{2\pi i/3},e^{2\pi i/3},e^{-2\pi i/3},e^{-2\pi i/3},e^{2\pi i/3}),\\
G(\kappa)&=&\mbox{diag}\,(J,J,J,J)
\end{array}
\end{equation}
where $J=\left[\begin{smallmatrix} 0 & 1\\1 & 0\end{smallmatrix}\right]$.
From equation~(\ref{c_act}), we know that $\Psi$ commutes with
the actions~(\ref{act3}) and~(\ref{act8}). Now, recall that $\Psi$ is 
chosen such that $(\Psi,\Phi)=I$ and note that $\phi_3=\ov{\phi}_1$, 
$\phi_4=\ov{\phi}_2$, $\phi_7=\ov{\phi}_5$, $\phi_8=\ov{\phi}_6$, then 
along with the equivariance condition on $\Psi$ we obtain
\[
\Psi(0)=
\left(\begin{array}{c}
  \mu (1  \quad \ov{\omega} \quad  \omega)\\
  \mu (1  \quad \omega \quad  \ov{\omega})\\
  \ov{\mu} (1  \quad \omega  \quad \ov{\omega})\\
  \ov{\mu} (1  \quad \ov{\omega}  \quad \omega)\\
  \zeta (1  \quad \ov{\omega}  \quad \omega)\\
  \zeta (1  \quad \omega  \quad \ov{\omega})\\
  \ov{\zeta} (1  \quad \omega  \quad  \ov{\omega})\\
  \ov{\zeta} (1  \quad  \ov{\omega} \quad \omega)
\end{array}\right),
\]
where $\mu$ and $\zeta$ are some nonzero complex numbers.

\paragraph{Basis of the unfolding space.}
Using~(\ref{dimW}), the dimension of the unfolding space 
${\cal W}$ for matrix $B$ is $16$ since here $r=4$ is the number of 
distinct eigenvalues, $k_{j}=2$ is the number of Jordan blocks for 
each eigenvalue and $n_{j,\ell}=1$ is the size of each Jordan block. 

We now define the $16$ matrices which form the basis of ${\cal E}({\cal W})$.
We define $3\times 8$ matrices $R_{j;\xi,\vartheta,m}$ for each 
eigenvalue $\lambda_j$, where $\lambda_1=\omega_1 i$, $\lambda_2=-\omega_1 i$,
$\lambda_3=\omega_2 i$, $\lambda_4=-\omega_2 i$. 
Let $v_i=(v_i^1,v_i^2,v_i^3)^t$, $i=1,\ldots,8$ be $3\times 1$ complex 
vectors. The linear mappings $\Pi_{j}$ are
\begin{equation}\label{proj}
\begin{array}{cc}
\Pi_1(v)=\left(\begin{array}{c}
\mu(v^1+\ov{\omega}v^2+\omega v^3)\\
\mu(v^1+\omega v^2+\ov{\omega} v^3)
\end{array}\right) & 
\Pi_2(v)=\left(\begin{array}{c}
\ov{\mu}(v^1+\omega v^2+\ov{\omega} v^3)\\
\ov{\mu}(v^1+\ov{\omega} v^2+\omega v^3)
\end{array}\right)\\
\\
\Pi_3(v)=\left(\begin{array}{c}
\zeta(v^1+\ov{\omega}v^2+\omega v^3)\\
\zeta(v^1+\omega v^2+\ov{\omega} v^3)
\end{array}\right) & 
\Pi_4(v)=\left(\begin{array}{c}
\ov{\zeta}(v^1+\omega v^2+\ov{\omega} v^3)\\
\ov{\zeta}(v^1+\ov{\omega}v^2+\omega v^3)
\end{array}\right).
\end{array}
\end{equation}

We choose $\xi,\vartheta\in\{1,2\}$ and $m=1$ since each Jordan
block is one-dimensional. Then, $R_{j;\xi,\vartheta,m}$ has vector 
$v_{2(j-1)+\xi}$ in column $2(j-1)+\vartheta$ and all the other columns 
are zero where
\begin{equation}\label{v-eqs}
\Pi_{j}(v_{2j-1})=\left(\begin{array}{c}
  1 \\ 0
\end{array}\right), \quad
\Pi_j(v_{2j})=\left(\begin{array}{c}
  0 \\ 1
\end{array}\right)
\end{equation}
and $v_3=\ov{v_1}$, $v_4=\ov{v_2}$, $v_7=\ov{v_5}$, $v_8=\ov{v_6}$.

To simplify the notation, we order the $16$ matrices $R_{j;\xi,\vartheta,m}$
in the following way:
\[
\{R^{1+8(k-1)},\ldots,R^{8+8(k-1)}\}=
\{R_{1;1,k,1},R_{2;1,k,1},R_{3;1,k,1},R_{4;1,k,1},
R_{1;2,k,1},R_{2;2,k,1},R_{3;2,k,1},R_{4;2,k,1}\}
\]
for $k=1,2$. One can verify using a computer algebra package that
\[
\mbox{rank}(\Phi(\tau_0),\Phi(\tau_1),\Phi(\tau_2),\Phi(\tau_3))=8,
\]
therefore we need four distinct delays $\tau_0,\tau_1,\tau_2,\tau_3$ 
to solve
\begin{equation}\label{RA-eqs}
R^{m}=\dps\sum_{i=0}^{3} A_{i}^{m} \Phi(\tau_{i}).
\end{equation}

\paragraph{Projection of the basis.}
We project the matrices $R^{m}$ on the subspace of $\D_3$-equivariant 
$3\times 8$ matrices, using the actions on $\C^{3}$~(\ref{act3}) and 
$\C^{8}$~(\ref{act8}). 
\begin{prop}
Let $\nu=(\nu^1,\nu^2,\nu^3)^{t}\in \C^{3}$ and consider the $3\times 8$ 
matrix $N_p$ having $\nu$ as its $p^{th}$ column and all the other columns 
are zero. From~(\ref{newproj}), we obtain a projection operator from which 
we compute
\[
\ov{N}_p=\sum_{g \in\D_3} \rho(g) N_p G(g^{-1}).
\]
Moreover, let $\eta_1(\nu)=\nu^1+\ov{\omega}\nu^2+\omega \nu^3$, 
$\eta_2(\nu)=\nu^1+\omega \nu^2+\ov{\omega} \nu^3$ and 
recall that $v=(1,\omega,\ov{\omega})^{t}$.
\begin{enumerate}
\item If $p=1,5$, then $\ov{N}_{p}$ has only nonzero columns $p$ and $p+1$ 
given respectively by $\eta_1 v$ and $\eta_1 \ov{v}$,

\item if $p=4,8$, then $\ov{N}_{p}$  has only nonzero columns $p-1$ and $p$ 
given respectively by $\eta_1\ov{v}$ and $\eta_1 v$,

\item if $p=3,7$, then $\ov{N}_{p}$  has only nonzero columns $p$ and $p+1$ 
given respectively by $\eta_2 \ov{v} $ and $\eta_2 v$,

\item and if $p=2,6$, then $\ov{N}_{p}$  has only nonzero columns $p-1$ 
and $p$ given respectively by $\eta_2 v$ and 
$\eta_2 \ov{v}$.
\end{enumerate}
\label{prop:R}
\end{prop}

\proof The proof is a straightforward calculation. \qed

From the definition of the matrices $R^{m}$, for $m=1,\ldots,8$ we have
$\vartheta=1$ and for $m=9,\ldots,16$ we have $\vartheta=2$. 
Table~\ref{tab:RN} shows the correspondence of matrices $R^{m}$ with matrices
$N_p$ of Proposition~\ref{prop:R} as well as the information given
by~(\ref{proj}) and~(\ref{v-eqs}) needed to compute $\ov{R}^{m}$.
From Proposition~\ref{prop:R} and Table~\ref{tab:RN}, we obtain easily that 
$\ov{R}^{m}=0$ for $m=5,\ldots,12$ and
\[
\begin{array}{cc}
\ov{R}^1=\ov{R}^{13}=[\mu^{-1}v,\mu^{-1}\ov{v},0,0,0,0,0,0], &
\ov{R}^2=\ov{R}^{14}=[0,0,\ov{\mu}^{-1}\ov{v},\ov{\mu}^{-1}v,0,0,0,0],\\
\ov{R}^3=\ov{R}^{15}=[0,0,0,0,\zeta^{-1}v,\zeta^{-1}\ov{v},0,0], &
\ov{R}^4=\ov{R}^{16}=[0,0,0,0,0,0,\ov{\zeta}^{-1}\ov{v},\ov{\zeta}^{-1}v].
\end{array}
\]

\arraystart
\begin{table}[htb]
\begin{center}
\begin{tabular}{|c|c|c|c|c|}\hline
Correspondence $R^{m}$ with $N_p$ & Projection & $v_j$ & $\eta_1(v_j)$ &
$\eta_2(v_j)$\\\hline
$R^1\leftrightarrow N_1$ and $R^9\leftrightarrow N_2$ & $\Pi_1$ & $v_1$ &
$\mu^{-1}$ & 0\\
$R^2\leftrightarrow N_3$ and $R^{10}\leftrightarrow N_4$ & $\Pi_2$ & $v_3$ &
$0$ & $\ov{\mu}^{-1}$\\
$R^3\leftrightarrow N_5$ and $R^{11}\leftrightarrow N_6$ & $\Pi_3$ & $v_5$ &
$\zeta^{-1}$ & 0\\
$R^4\leftrightarrow N_7$ and $R^{12}\leftrightarrow N_8$ & $\Pi_4$ & $v_7$ &
$0$ & $\ov{\zeta}^{-1}$\\
$R^5\leftrightarrow N_1$ and $R^{13}\leftrightarrow N_2$ & $\Pi_1$ & $v_2$ &
$0$ & $\mu^{-1}$\\
$R^6\leftrightarrow N_3$ and $R^{14}\leftrightarrow N_4$ & $\Pi_2$ & $v_4$ &
$\ov{\mu}^{-1}$ & $0$\\
$R^7\leftrightarrow N_5$ and $R^{15}\leftrightarrow N_6$ & $\Pi_3$ & $v_6$ &
$0$ & $\zeta^{-1}$\\
$R^8\leftrightarrow N_7$ and $R^{16}\leftrightarrow N_8$ & $\Pi_4$ & $v_8$ &
$\ov{\zeta}^{-1}$ & 0\\\hline
\end{tabular}
\end{center}
\caption{The first column shows the correspondence between $R^{m}$ and
$N_p$ of Proposition~\ref{prop:R} while the remaining columns detail the
information needed to compute matrices $\ov{R}^{m}$ of the first column
where $\eta_j(\cdot)$ is defined in Proposition~\ref{prop:R}.}
\label{tab:RN}
\end{table}
\arrayfinish

\paragraph{$\D_3$-mini-versal unfolding}
Since ${\cal E}({\cal W})$ is an unfolding space for $B$, we again set 
$\hat{B}_{j}^{\D_3}=\Psi(0)\ov{R}^{j}$ for $j=1,\ldots,16$. 
The space ${\cal W}$ is spanned by sixteen $8\times 8$ matrices. For 
$i=1,\ldots 8$, $\Omega_i$ has a $1$ at element $(i,i)$ and zeroes elsewhere.
For $i=9,10,11,12$, $\Omega_i$ has a $1$ at element $(i-8,13-i)$ and
zeroes elsewhere while for $i=13,14,15,16$, $\Omega_i$ has a $1$ at
element $(i-8,21-i)$ and zeroes elsewhere. Now, obviously,
$\hat{B}_{j}^{\D_3}=0$ for $j=5,\ldots,12$ and a computation shows that
\[
\begin{array}{ll}
\hat{B}_1^{\D_3}=\hat{B}_{13}^{\D_3}=&[B,y_1]+(1+(1+\mu^{-1})|\omega|^2)\Omega_1+(1+2|\omega|^2)\Omega_2+\mu^{-1}\ov{\mu}(1+2|\omega|^2)(\Omega_{11}+\Omega_{12})
\\
\hat{B}_2^{\D_3}=\hat{B}_{14}^{\D_3}=&[B,y_2]+(1+2|\omega|^2)\Omega_3+(1+2|\omega|^2)\Omega_4+(\mu\ov{\mu}^{-1}(1+|\omega|^2)+\ov{\mu}^{-1}|\omega|^2)\Omega_9\\
&+\mu\ov{\mu}^{-1}(1+2|\omega|^2)\Omega_{10}
\\
\hat{B}_3^{\D_3}=\hat{B}_{15}^{\D_3}=&[B,y_3]+(1+(1+\zeta^{-1})|\omega|^2)\Omega_5+(1+2|\omega|^2)\Omega_6+\zeta^{-1}\ov{\zeta}(1+2|\omega|^2)(\Omega_{15}+
\Omega_{16})
\\
\hat{B}_4^{\D_3}=\hat{B}_{16}^{\D_3}=&[B,y_4]+(1+2|\omega|^2)\Omega_7+(1+2|\omega|^2)\Omega_8+(\zeta\ov{\zeta}^{-1}(1+|\omega|^2)+\ov{\zeta}^{-1}|\omega|^2)\Omega_{13}\\
&+\zeta\ov{\zeta}^{-1}(1+2|\omega|^2)\Omega_{14}.
\end{array}
\]
for some matrices $y_1,y_2,y_3,y_4\in \mbox{Mat}_{8\times 8}$. 
The four first rows of $\Theta$ are 
\[
\begin{array}{l}
v_1=(1+(1+\mu^{-1})|\omega|^2,1+2|\omega|^2,0,0,0,0,0,0,0,0,\mu^{-1}\ov{\mu}(1+2|\omega|^2),\mu^{-1}\ov{\mu}(1+2|\omega|^2),0,0,0,0)\\
v_2=(0,0,1+2|\omega|^2,1+2|\omega|^2,0,0,0,0,\ov{\mu}^{-1}(\mu(1+|\omega|^2)+|\omega|^2),\mu\ov{\mu}^{-1}(1+2|\omega|^2),0,0,0,0,0,0)\\
v_3=(0,0,0,0,1+(1+\zeta^{-1})|\omega|^2,1+2|\omega|^2,0,0,0,0,0,0,0,0,\zeta^{-1}\ov{\zeta}(1+2|\omega|^2),\zeta^{-1}\ov{\zeta}(1+2|\omega|^2))\\
v_4=(0,0,0,0,0,0,1+2|\omega|^2,1+2|\omega|^2,0,0,0,0,\ov{\zeta}^{-1}(\zeta(1+|\omega|^2)+|\omega|^2),\zeta\ov{\zeta}^{-1}(1+2|\omega|^2),0,0),
\end{array}
\]
rows $5$ to $12$ are zero and rows $13$ to $16$ are identical to rows
$1$ to $4$. Thus, $\{v_1,v_2,v_3,v_4\}$ determine a maximal set 
of linearly independent vectors in the rowspace of $\Theta$. Hence,
\[
\mbox{Mat}_{8\times 8}^{\D_3}=T_{B}\Sigma^{\D_3}\oplus 
\C\cdot\{\hat{B}_1^{\D_3},\hat{B}_2^{\D_3},\hat{B}_3^{\D_3},
\hat{B}_4^{\D_3}\}
\]
and we need only consider matrices $\ov{R}^{j}$ for $j=1,2,3,4$
in order to compute the $\D_3$-equivariant $\Lambda$-mini-versal unfolding.

\paragraph{Computation of the matrices $\mathbf{A}_j^m$.} 
We take $\D_3$-equivariant
matrices $A_{j}^m$ as in~(\ref{D3mat3}) and we compute
\[
A_{j}^{m}\Phi(\tau_{j})=(a_j^m+b_j^m(\ov{\omega}+\omega))\Phi(\tau_j).
\]
Thus, equations $\ov{R}^{m}=\sum_{j=0}^{3} A_{j}^{m} \Phi(\tau_{j})$
for $m=1,\ldots,4$ can be solved by finding the solution
to the first row only.
\begin{equation}\label{2eqs}
\left(\begin{array}{c}
\mu^{-1}\\ 0 \\ 0 \\ 0
\end{array}\right)
=M 
\left(\begin{array}{c}
a_0^1+b_0^1(\ov{\omega}+\omega)\\ a_1^1+b_1^1(\ov{\omega}+\omega) \\ 
a_2^1+b_2^1(\ov{\omega}+\omega) \\ a_3^1+b_3^1(\ov{\omega}+\omega)
\end{array}\right)
\end{equation}
where $M$ is the generically nonsingular matrix~(\ref{matM}).

\begin{rmk}\label{remB}
Note that, as in the previous example, from the form of 
equations~(\ref{2eqs}) we can choose $a_{j}^{2}=\ov{a}_{j}^1$, 
$a_{j}^{4}=\ov{a}_{j}^{3}$, $b_{j}^{2}=\ov{b}_{j}^{1}$ and 
$b_{j}^{4}=\ov{b}_{j}^{3}$ for $j=0,1,2,3$.
\end{rmk}

\paragraph{$\D_3$-equivariant $\Lambda$-unfolding.}
For $m=1,2,3,4$, define
\[
L_{m}\,z=\sum_{j=0}^{3} A_{j}^{m}z(\tau_j)=\sum_{j=0}^{3}
\left(
\begin{array}{ccc}
a_{j}^{m} & b_{j}^{m} & b_{j}^{m} \\
   b_{j}^{m} & a_{j}^{m} & b_{j}^{m}\\
   b_{j}^{m} & b_{j}^{m} & a_{j}^{m}
\end{array}
\right)z(\tau_j),
\]
and then the complex $\D_3$-equivariant $\Lambda$-versal unfolding is defined by:
\[
\begin{array}{rcl}
{\cal L}(\alpha)\,z&=&{\cal L}_{0}\,z
+\dps\sum_{j=0}^{3}\left( \dps\sum_{m=1}^{4} \alpha_{m}
\left(\begin{array}{ccc}
a_{j}^{m} & b_{j}^{m} & b_{j}^{m}\\
b_{j}^{m} & a_{j}^{m} & b_{j}^{m}\\
b_{j}^{m} & b_{j}^{m} & a_{j}^{m}
\end{array}\right)\right)z(\tau_j).
\end{array}
\]
As in the previous example, setting $\tau_0=0$, $\tau_1=\tau_s^{*}$, 
$\tau_2=\tau_n^{*}$ and $b_0^m=b_1^m=a_2^m=a_3^m=0$ for $m=1,2,3,4$ 
preserves the structure of the delay-differential equation.
Then as a delay-differential equation, the complex unfolding is given 
exactly by~(\ref{unfDDE}) with again
$\epsilon_1=\sum_{m=1}^{4} \alpha_m a_{0}^{m}$,
$\epsilon_2=\sum_{m=1}^{4} \alpha_m a_{1}^{m}$,
$\epsilon_3=\sum_{m=1}^{4} \alpha_m b_{2}^{m}$, and
$\epsilon_4=\sum_{m=1}^{4} \alpha_m b_{3}^{m}$.

Note that the unfolding parameters $\epsilon_1,\epsilon_2,\epsilon_3,
\epsilon_4$ are linearly independent by Remark~\ref{remB} and from
the same calculation as in the previous example. Moreover, from 
Remark~\ref{remB}, the real $\D_3$-equivariant $\Lambda$-versal unfolding 
is also given by~(\ref{unfDDE}).
\vspace*{0.2in}

\noindent
{\Large\bf Appendix}

\appendix
\Section{Proof of Lemma~\ref{suff_cond_versal_prop_gamma}}

We will prove the result in the case where (\ref{suff_cond_versal_gamma}) holds 
with $p$ equal to the codimension of $T_B\Sigma^{\Gamma}$ in 
$\mbox{\rm Mat}^{\Gamma}_{c\times c}$.
The proof follows along the same lines as the non-equivariant case (which can be
found, for example in \cite{Ar}).  
We first define the $\Gamma$-centralizer of
$u\in\mbox{\rm Mat}^{\Gamma}_{c\times c}$ as the linear subspace $Z^{\Gamma}_u$ of 
all matrices in $\mbox{\rm Mat}^{\Gamma}_{c\times c}$ which commute with $u$.
Consider now the $\Gamma$-centralizer $Z^{\Gamma}_{B}$ of the matrix $B$.  In the
space $GL(c,{\mathbb C})\cap \mbox{\rm Mat}^{\Gamma}_{c\times c}$, let ${\cal P}$ 
be a smooth surface transversal to $I_{c\times c}+Z^{\Gamma}_{B}$ at $I_{c\times c}$,
and such that the dimension of ${\cal P}$ is equal to the codimension 
of $Z^{\Gamma}_{B}$ in $\mbox{\rm Mat}^{\Gamma}_{c\times c}$.  We then define a mapping
\[
{\cal F}:{\cal P}\times {\mathbb C}^p\longrightarrow \mbox{\rm Mat}^{\Gamma}_{c\times c}\,\,,
\,\,\,\,\,\,\,\,\,\,\,\,\,\,\,\,\,
{\cal F}(p,\alpha)=p\,{\cal B}(\alpha)\,p^{-1}.
\]
The key result we need is
\begin{lemma}
${\cal F}$ is a local diffeomorphism in a neighborhood of $(I_{c\times c}\,,\,\alpha_0)$.
\label{keylem}
\end{lemma}
\proof
First, we note that the mapping
\[\begin{array}{c}
f:GL(c,{\mathbb C})\cap \mbox{\rm Mat}^{\Gamma}_{c\times c}\longrightarrow
\Sigma^{\Gamma}\subset\mbox{\rm Mat}^{\Gamma}_{c\times c}\,\,,\\[0.15in]
f(C)=CBC^{-1}
\end{array}
\]
is such that the derivative of $f$ evaluated at $I_{c\times c}$ is
\[
\begin{array}{c}
Df(I_{c\times c}):T_{I_{c\times c}}\mbox{\rm Mat}^{\Gamma}_{c\times c}\cong
\mbox{\rm Mat}^{\Gamma}_{c\times c}\longrightarrow 
T_B\Sigma^{\Gamma}\subset
T_B\mbox{\rm Mat}^{\Gamma}_{c\times c}
\cong\mbox{\rm Mat}^{\Gamma}_{c\times c}
\,\,,\\[0.15in]
Df(I_{c\times c})(u)=[u,B],
\end{array}
\]
which is onto $T_B\Sigma^{\Gamma}$.
It follows from the above that the dimension of $Z^{\Gamma}_B$ is equal to the
codimension of $T_B\Sigma^{\Gamma}$ in $\mbox{\rm Mat}^{\Gamma}_{c\times c}$, which is
equal to $p$.  Also, the dimension of ${\cal P}$ 
is equal to the dimension of $T_B\Sigma^{\Gamma}$.

Now, the mapping ${\cal F}$ is such that
\[
D_p{\cal F}(I_{c\times c},\alpha_0)(u,\alpha)=Df(I_{c\times c})(u)=[u,B]\,\,,\,\,\,\,\,\,\,\,\,\,\,\,\,\,\,\,\,
D_{\alpha}{\cal F}(I_{c\times c},\alpha_0)(u,\alpha)=D_{\alpha}{\cal B}(\alpha_0)\alpha\,\,.
\]
By construction of ${\cal P}$, we have that
\[
D_p{\cal F}(I_{c\times c},\alpha_0)=\left.Df(I_{c\times
    c})\right|_{T_{I_{c\times c}}{\cal P}}
\] 
maps $T_{I_{c\times c}}{\cal P}$ isomorphically onto $T_B\Sigma^{\Gamma}$.
Also, by the hypothesis of 
Lemma~\ref{suff_cond_versal_prop_gamma}, $D_{\alpha}{\cal
  B}(\alpha_0)$
maps $T_{\alpha_0}{\mathbb C}^p\cong {\mathbb C}^p$ isomorphically
onto a space which is a direct sum complement of $T_B\Sigma^{\Gamma}$
in $\mbox{\rm Mat}^{\Gamma}_{c\times c}$.  It thus follows that
$D{\cal F}(I_{c\times c},\alpha_0)$ is an isomorphism between vector
spaces of dimension $\mbox{\rm dim}(\mbox{\rm Mat}^{\Gamma}_{c\times c})$, 
so ${\cal F}$ is a local
diffeomorphism by the inverse function theorem. 
\hfill
\qed

Now, let ${\cal A}(\beta)$ be a $q$-parameter $\Gamma$-unfolding
of $B$ (with ${\cal A}(\beta_0)=B$).  
Define $\Pi_1$ and $\Pi_2$ as the projections of ${\cal P}\times
{\mathbb C}^p$ onto ${\cal P}$ and ${\mathbb C}^p$ respectively.

For all $\beta$ sufficiently close to $\beta_0$ in ${\mathbb C}^q$, define
\[\begin{array}{c}
C(\beta)=\Pi_1\,{\cal F}^{-1}({\cal A}(\beta))\\[0.15in]
\phi(\beta)=\Pi_2\,{\cal F}^{-1}({\cal A}(\beta)).
\end{array}
\]
It follows that for all $\beta$ sufficiently close to $\beta_0$, 
\[
{\cal A}(\beta)={\cal F}(C(\beta),\phi(\beta))=C(\beta)\,{\cal B}(\phi(\beta))
\,(C(\beta))^{-1},
\]
which proves the lemma.

\vspace*{0.25in}
\noindent
{\Large\bf Acknowledgements}

\vspace*{0.2in}
This research is partly supported by the
Natural Sciences and Engineering Research Council of Canada in the
form of a postdoctoral fellowship and the CRM (PLB) and an 
individual research grant (VGL).  The authors are grateful to David 
Handelman for some very helpful suggestions.

\end{document}